\documentclass[twoside,a4paper,10pt]{article}
\usepackage{latexsym,amssymb,amsmath}
\usepackage{hyperref}
\usepackage{graphicx}
\usepackage{enumerate}
\usepackage[frenchb,english]{babel}
\usepackage[matrix,arrow]{xy}
\usepackage{color}
\def \be{\begin{eqnarray*}}
\def \ee{\end{eqnarray*}}
\def \ben{\begin{enumerate}}
\def \een{\end{enumerate}}
\def \beit{\begin{itemize}}
\def \eeit{\end{itemize}}
\def \bui#1#2{\mathrel{\mathop{\kern 0pt#1}\limits^{#2}}}
\def \buil#1#2{\mathrel{\mathop{\kern 0pt#1}\limits_{#2}}}
\def \bfll{\begin{flushleft}}
\def \efll{\end{flushleft}}
\def \bflr{\begin{flushright}}
\def \eflr{\end{flushright}}

\def \findemo{\hfill$\square$\\}

\def \ovl{\overline}

\def \wit{\widetilde}
\def \wnabla{\wit{\nabla}}

\def \bg{\overline{g}}

\def \C{\mathbb{C}}

\def \R{\mathbb{R}}

\newtheorem{ethm}{Theorem}[section]

\newtheorem{elemme}[ethm]{Lemma}
\newtheorem{erem}[ethm]{Remark}

\newtheorem{ecor}[ethm]{Corollary}
\newtheorem{prop}[ethm]{Proposition}
\newtheorem{eexemple}[ethm]{Example}

\title{An Obata-type characterization of doubly-warped product K\"ahler manifolds}
\author{Nicolas Ginoux\footnote{Universit\'e de Lorraine, CNRS, IECL, F-57000 Metz, France, \texttt{nicolas.ginoux@univ-lorraine.fr}
}, Georges Habib\footnote{Lebanese University, Department of Mathematics, P.O. Box 90656 Fanar-Matn, Lebanon, \texttt{ghabib@ul.edu.lb}}, Mihaela Pilca\footnote{Universit\"at Regensburg, Universit\"atstr. 31, 93051 Regensburg, Germany, \texttt{mihaela.pilca@ur.de}} and Uwe Semmelmann\footnote{
Universit{\"a}t Stuttgart,
Pfaffenwaldring 57,
70569 Stuttgart, Germany, \texttt{uwe.semmelmann@mathematik.uni-stuttgart.de}}}

\begin{document}
\maketitle

\noindent\begin{center}\begin{tabular}{p{155mm}}
\begin{small}{\bf Abstract.} 
We give a characterization {\sl \`a la Obata} for certain families of K\"ahler manifolds. These results are in the same line as other extensions of the well-known Obata's rigidity theorem from \cite{Obata62}, like for instance the generalizations in \cite{RanjSant97} and \cite{Santhanam07}. Moreover, we give a complete description of the so-called K\"ahler doubly-warped product structures whose underlying metric is Einstein.
\end{small}\\
\end{tabular}\end{center}

\section{Introduction}\label{s:intromainresults}

This paper is the first of two papers devoted to the classification of complete K\"ahler manifolds carrying a real-valued function $u$ whose Hessian  is $J$-invariant and has pointwise at most two eigenvalues and one of them has as eigenvector the gradient of $u$. In this first paper we consider the case where $u$ has {\it no critical point}, and in the forthcoming paper \cite{GinHabPilSemm19} we treat the critical case, for which further constructions are needed. Before stating our main result, let us review some previous results that motivate our study. \\

In \cite[Theorem p. 614]{Obata62}, it is shown that the only complete Riemannian manifold $(M^n,g)$ carrying a real-valued function $u$ whose Hessian satisfies $\nabla^2 u=-u {\rm Id}$ is the round sphere. This result, known as the Obata theorem, has been generalized on K\"ahler manifolds in several papers such as \cite{MolzonPinneyMortensen93, RanjSant97, Santhanam07}. Namely, in \cite[Theorem 3]{RanjSant97}, the authors proved that a complete K\"ahler manifold $(M^{2n},g,J)$ is biholomorphically isometric to $\C\mathrm{P}^n$ with holomorphic sectional curvature $1$ if and only if there exists a function $u$ whose Hessian has at most two eigenvalues, namely $-\frac{u+1}{2}$ and $-u$ and where $\nabla u$ is an associated eigenvector (see also \cite[Theorem p. 614]{MolzonPinneyMortensen93} for a weaker version). In \cite[Theorem 1]{Santhanam07}, G. Santhanam proved that given a function $u$ on a complete K\"ahler manifold $(M^{2n},g,J)$  whose Hessian has the eigenvalues $u$ and $\frac{u+1}{2}$ and where $\nabla u$ and $J\nabla u$ are both eigenvectors associated to $u$, then the manifold $M$  is either isometric to the complex hyperbolic space $\C\mathrm{H}^n$ of constant sectional curvature $-1$ or it is diffeomorphic to the normal bundle of some $2$-codimensional totally geodesic submanifold $M_0$ of $M$, such that the fibre over each point in $M_0$ is isometric to the hyperbolic plane $\mathbb{H}^2$ of constant curvature $-1$.\\   

The main result of this paper is the following:  

\begin{ethm}\label{t:cardwp2eigenv}
Let $(\wit{M}^{2n},\wit{g},\wit{J})$ be a connected complete K\"ahler manifold of real dimension $2n\geq4$ carrying a function $u\in C^\infty(\wit{M},\mathbb{R})$ without critical points which satisfies the following two conditions:
\begin{enumerate}[$\bullet$]
\item its Hessian $\wnabla^2 u$ is $\wit{J}$-invariant and has pointwise at most two eigenvalues $\lambda$ and $\mu$,
\item its gradient $\wnabla u$ is a pointwise eigenvector of $\wnabla^2u$ with the eigenvalue $\lambda$.
\end{enumerate}
Then the following claims hold true:
\begin{enumerate}[i)]
\item If $\mu$ vanishes at one point of $\wit M$, then $\mu$ vanishes identically on $\wit M$ and the triple $(\wit{M}^{2n},\wit{g},\wit{J})$ is locally biholomorphically isometric to $(\R_t\times\R_s\times\Sigma,dt^2\oplus\rho^2(t)ds^2\oplus g_{\Sigma})$ for some K\"ahler manifold $(\Sigma^{2n-2},g_{\Sigma})$ and  $\rho(t):=|\wnabla u|(t,s,x)$, where the complex structure of $(\R_t\times\R_s\times\Sigma,dt^2\oplus\rho(t)^2ds^2\oplus g_{\Sigma})$ is the one that is naturally induced by the complex structure of $(\Sigma^{2n-2},g_{\Sigma})$.
\item If $\mu$ does not vanish at any point of $\wit M$, then we distinguish the following two cases:
\begin{enumerate}
\item If $n>2$, then up to changing $u$ into $a u+b$ with $a,b\in\R$, $a\neq 0$, the function $u$ may be assumed to be positive and the K\"ahler manifold $(\wit{M}^{2n},\wit{g},\wit{J})$ is biholomorphically isometric to a doubly-warped product $\left(\R\times M^{2n-1},dt^2\oplus\rho(t)^2\left(\rho'(t)^2\hat{g}_{\hat{\xi}}\oplus\hat{g}_{\hat{\xi}^\perp}\right)\right)$, where $M$ is a level hypersurface of $u$, the triple $(M^{2n-1},\hat{g},\hat{\xi})$ is Sasaki and $\rho(t)=\sqrt{u(t,x)}$, for any $(t,x)\in\R\times M$.
\item If $n=2$, then
up to changing $u$ into $au+c$ with $a,b\in\R$, $a\neq 0$,  the function $u$ must be positive and the K\"ahler manifold $(\wit{M}^{2n},\wit{g},\wit{J})$ is biholomorphically isometric to a doubly-warped product $\left(\R\times M^{2n-1},dt^2\oplus\rho(t)^2\left(\rho'(t)^2\hat{g}_{\hat{\xi}}\oplus\hat{g}_{\hat{\xi}^\perp}\right)\right)$, where $M$ is a level hypersurface of $u$, the triple $(M^{2n-1},\hat{g},\hat{\xi})$ is a minimal Riemannian flow that is basic conformally Sasaki and $\rho(t)=\sqrt{u(t,x)}$, for any $(t,x)\in\R\times M$.
\end{enumerate}
Moreover, in this case ($\mu\neq0$), we have that $\lambda\circ F(t,x)=\frac{\partial^2 (u\circ F)}{\partial t^2}(t,x)$ and $\mu=\frac{|\wnabla u|^2}{2u}$.
\end{enumerate}
\end{ethm}

The assumptions of Theorem  \ref{t:cardwp2eigenv}  are related to various other well studied
situations. First, it is easy to check that the condition of a 
$\wit{J}$-invariant Hessian 
$\tilde \nabla^2 u$ is equivalent to the condition that $K:= - \wit{J} \tilde 
\nabla u$
is a Hamiltonian Killing vector field with moment map $u$, i.e. we  have
$L_K\wit{J}=0 = L_K \tilde g$ and $K \lrcorner \, \omega = du$, where $\omega$ 
denotes the K\"ahler form of $(\wit{M},\wit{g},\wit{J})$.

\medskip

Next, the condition that the gradient $\tilde \nabla u$ is a pointwise eigenvector of the Hessian 
of $u$, say for some eigenvalue $\lambda$, is equivalent to the equation
$dx = 2 \lambda du $, where $x$ is the length function $x=|K|^2$. In particular,
$dx \wedge du = 0$ and $x$ has to be a function of $u$. Then the local $S^1$-action 
generated by $K$ is {\it rigid} in the sense of V. Apostolov et al. (cf. \cite{ACG}).

\medskip

Our main theorem is also related to the work of A. Derdzinski and G. Maschler in 
\cite{DM}, where they studied the question whether a given K\"ahler metric is 
conformal  to an Einstein metric. As a necessary condition  for the conformal 
factor $u$ they obtained that $\wit{J} \wnabla u$ has to be a Killing vector 
field and an eigenvector of the Hessian of
$u$ and of  the Ricci tensor.
They called such functions special K\"ahler-Ricci potentials.

\medskip

Another equivalent way of formulating the assumptions of Theorem  
\ref{t:cardwp2eigenv} is in terms of
the distribution $\mathcal D_+$ spanned by $K$ and $JK$. It turns out that equivalently  this distribution has 
to be totally geodesic, holomorphic and conformal. The last condition follows 
from
the assumption that the Hessian of $u$ has at most two eigenvalues. Then our 
metric $\tilde g$ is locally of Calabi type and there is a local classification 
due to  S. Chiossi and P.-A. Nagy in \cite{ChiossiNagy10}. 
From this point of view  it becomes
clear that our  manifolds are also {\it ambi-K\"ahler}, i.e. switching the sign
of the complex structure $\wit{J}$ along the distribution $\mathcal D_+$ 
defines a new
integrable complex structure $I$ such that $(u^{-2}g, I)$ is again K\"ahler.
We recall that there
is a classification of compact ambi-K\"ahler manifolds in the work of F. Madani, A. Moroianu
and M. Pilca in \cite{MMP}. 
Moreover, $u \omega_+$ is a Hamiltonian $2$-form of rank $1$, where $\omega_+$
is the restriction of the K\"ahler form to $\mathcal D_+$.  
Manifolds admitting Hamiltonian $2$-forms are studied in a series of papers of V. Apostolov et al. 
including a global classification in the compact case (cf. \cite{ACG}, \cite{ACGT}).

\medskip

In contrast to the results mentioned so far our main theorem gives a global
description of the manifold without the additional compactness assumption. 
In fact, as a result, the underlying manifold in Theorem \ref{t:cardwp2eigenv} cannot be compact.\\

The idea of the proof of Theorem \ref{t:cardwp2eigenv} consists in identifying the manifold $\wit{M}$ with $I\times M$,  where $M$ denotes a level hypersurface of $u$,   via the flow of the normalized gradient $\nu:=\frac{\wnabla u}{|\wnabla u|}\in\Gamma(T\wit{M})$, which is geodesic. 
We show that the vector field $\xi:=-\wit{J}\nu$ defines a minimal Riemannian flow on $(M, \wit{g}\big|_{T^*M\otimes T^*M})$, whose O'Neill tensor coincides with the complex structure $\wit{J}$, up to some factor depending on the eigenvalue $\mu$.\\

Conversely, given any  K\"ahler doubly-warped product (see Lemma \ref{l:charactKaehlerrhosigmak} for the existence of such a structure) of the form $(I\times M^{2n-1},dt^2\oplus\rho^2((\rho')^2\widehat g_{\hat{\xi}}\oplus \widehat g_{\hat{\xi}^\perp}))$, where $I\subset \R,$  $\rho,\rho'\colon I\to\R$ are positive functions and $(M,\hat{g}=\widehat g_{\hat{\xi}}\oplus \widehat g_{\hat{\xi}^\perp},\hat{\xi})$ is Sasaki, a direct computation shows that the function $u:=\rho^2$ satisfies the following second-order PDE:
\[\wnabla^2u=\wnabla^2u(\nu,\nu)\cdot\left(\nu^\flat\otimes\nu+\xi^\flat\otimes\xi\right)+\frac{|\wnabla u|^2}{2u}\mathrm{Id}_{\{\xi,\nu\}^\perp},\]
where  $\wnabla u:=\mathrm{grad}_{\wit g}^{\wit{M}}(u)$, $\wnabla^2u:=\mathrm{Hess}_{\wit g}^{\wit{M}}(u)$, $\nu:=\frac{\wnabla u}{|\wnabla u|}$ and $\xi:=-\wit J\nu$. Hence the Hessian of $u$ has two eigenvalues that coincide with $\lambda$ and $\mu$ in Theorem \ref{t:cardwp2eigenv}. Note that the function $\rho$ itself has no $\wit J$-invariant Hessian, whereas $\rho^2$ does, that is why we consider $\rho^2$.\\ 

The paper is organized as follows. In Section \ref{ss:doublywarpedproducts}, we review some basic facts on doubly-warped products and characterize those which are K\"ahler. We recall that these structures were first introduced by Baier in his master thesis \cite{Baier97}, in order to compute the Dirac spectrum of the complex hyperbolic space. In Section \ref{s:Obatadwp}, we provide the proof of the main theorem. In the last section of the paper we investigate when the metric of a K\"ahler doubly-warped product is Einstein and discuss the solutions of the differential equation that the warping function $\rho$ has to satisfy. 

\section{K\"ahler doubly-warped products}\label{ss:doublywarpedproducts}

In this section, we recall some basic facts on doubly-warped products. We characterize among these manifolds those which are K\"ahler and provide the necessary integrability conditions. We refer to \cite{GinSemm11} for more details.\\

Let $M$ be a manifold and consider the product \cite[Lemma 3.1]{GinSemm11} 
\[(\wit{M}:=I\times M,\wit{g}:=\beta dt^2\oplus g_t),\]
where $I\subset\R$ is an open interval, $g_t$ is a smooth $1$-parameter family of Riemannian metrics on $M$ and $\beta\in C^\infty(I\times M,\R_+^\times)$. We can easily see that the Koszul formula implies the following identities for all $X,Y\in\Gamma(\pi_2^*TM)$, where $\pi_2: \wit{M}\to M$ \color{black} denotes the projection on the second factor: 
\begin{eqnarray}\label{equationwarpedproduct}
\wnabla_{\partial_t}\partial_t&=&-\frac{1}{2}\mathrm{grad}_{g_t}(\beta(t,\cdot))+\frac{1}{2\beta}\frac{\partial\beta}{\partial t}\partial_t,\\
\wnabla_{\partial_t}X&=&\frac{\partial X}{\partial t}+\frac{1}{2}g_t^{-1}\frac{\partial g_t}{\partial t}(X,\cdot)+\frac{1}{2\beta}\frac{\partial\beta}{\partial x}(X)\partial_t,\nonumber\\
\wnabla_X \partial_t&=&\frac{1}{2}g_t^{-1}\frac{\partial g_t}{\partial t}(X,\cdot)+\frac{1}{2\beta}\frac{\partial\beta}{\partial x}(X)\partial_t,\nonumber\\
\wnabla_X Y&=&\nabla_X^{M_t} Y-\frac{1}{2\beta}\frac{\partial g_t}{\partial t}(X,Y)\partial_t\nonumber,
\end{eqnarray}
where $\frac{\partial X}{\partial t}=[\partial_t,X]$ and $\nabla^{M_t}$ is the Levi-Civita covariant derivative of $(M,g_t)$. From now on, we assume $\hat{\xi}$ to be a unit Killing vector field, in other words $(M,\hat{g},\hat{\xi})$ is a so-called minimal Riemannian flow. For more details on Riemannian flows, we refer to \cite{Carriere84}. In this case, we have an orthogonal splitting $TM=\R\cdot\hat\xi\oplus\hat\xi^\perp$ and the normal bundle $Q:=\hat\xi^\perp$ of the flow admits a so-called transversal Levi-Civita connection, denoted by $\hat{\nabla}$, and which is defined for all  $X\in \Gamma(TM)$ and $Z\in\Gamma(Q)$ as follows, \emph{cf.} \cite{Tondeur88}:
$$\hat{\nabla}_XZ:=\left\{\begin{array}{ll}[\hat{\xi},Z]^Q,&\textrm{if }X=\hat{\xi}\\\left(\nabla_X^{\hat{M}}Z\right)^Q,&\textrm{if }X\in \Gamma(Q) ,\end{array}\right.$$ where $\left(\cdot\right)^Q$ denotes the $\hat{g}$-orthogonal projection $TM\to Q$ and $\nabla^{\hat{M}}$ denotes the Levi-Civita covariant derivative of $(M,\hat{g})$. The connection $\hat\nabla$ is compatible with the induced metric $\hat{g}_{\hat{\xi}^\perp}$ on the bundle $Q$ and it curvature vanishes along $\hat\xi$. Recall also that a minimal Riemannian flow is characterized by the fact that the tensor $\hat{h}:=\nabla^{\hat{M}}\hat{\xi}$, known as the O'Neill tensor \cite{Oneil66}, satisfies $\hat{h}(\hat{\xi})=0$ and is a skew-symmetric endomorphism field on $Q$ equal to $\hat{g}(\hat{h}(Y),Z)=-\frac{1}{2}\hat{g}([Y,Z],\hat{\xi})$ for any $Y,Z\in \Gamma(Q)$.\\

We consider in the following the general Ansatz
\begin{equation}\label{formofthemetric}
g_t:=\rho(t)^2\left(\sigma(t)^2\hat{g}_{\hat{\xi}}\oplus k(t,x)^2\hat{g}_{\hat{\xi}^\perp}\right),
\end{equation}
where $\rho,\sigma\colon I\to(0,\infty)$ and $k\colon I\times M\to(0,\infty)$ are {\sl a priori} arbitrary smooth positive functions and we define on the manifold $\wit{M}=I\times M$ the Riemannian metric $\wit{g}:=dt^2\oplus g_t$. In the next lemma, we explicit the Levi-Civita connection $\wnabla$ of $(\wit{M},\wit{g})$ and express it in terms of the transversal Levi-Civita connection~$\hat{\nabla}$.

\begin{elemme}\label{l:identitiesrhosigmak}
Let $(\wit{M},\wit{g}):=\left(I\times M,dt^2\oplus\rho(t)^2\left(\sigma(t)^2\hat{g}_{\hat{\xi}}\oplus k(t,x)^2\hat{g}_{\hat{\xi}^\perp}\right)\right)$, where $(M,\hat{g},\hat{\xi})$ is a minimal Riemannian flow.
Then for all $Z,Z'\in\Gamma(\pi_2^* Q)$, the following identities hold: 
$$\wnabla_{\partial_t}\partial_t=0, \quad 
\wnabla_{\partial_t}\xi=0, \quad
\wnabla_{\partial t}Z=\frac{\partial Z}{\partial t}+\frac{1}{\rho k}\cdot\frac{\partial(\rho k)}{\partial t}Z,$$
$$ \wnabla_\xi\partial_t=\frac{(\rho\sigma)'}{\rho\sigma}\xi, \quad
\wnabla_\xi\xi=-\frac{(\rho\sigma)'}{\rho\sigma}\partial_t, \quad
\wnabla_\xi Z=\hat{\nabla}_{\xi}Z+\frac{\xi(k)}{k}Z+\frac{\sigma}{\rho k^2}\hat{h}Z, $$
$$ \wnabla_Z\partial_t=\frac{1}{\rho k}\cdot\frac{\partial(\rho k)}{\partial t}Z, \quad \wnabla_Z\xi=\frac{\xi(k)}{k}Z+\frac{\sigma}{\rho k^2}\hat{h}Z \quad \text{and}$$
\begin{equation*}
\begin{split}
\wnabla_ZZ'=&\hat{\nabla}_ZZ'+Z(\ln(k))Z'+Z'(\ln(k))Z-\wit{g}(Z,Z')\mathrm{grad}^{\hat{\nabla}}(\ln(k))\\
&-\xi(\ln(k))\wit{g}(Z,Z')\xi-\frac{\sigma}{\rho k^2}\wit{g}(\hat{h}Z,Z')\xi-\frac{\partial\ln(\rho k)}{\partial t}\wit{g}(Z,Z')\partial_t,
\end{split}
\end{equation*}
where $\hat{h}:=\nabla^{\hat{M}}\hat{\xi}\in\Gamma(\mathrm{End}(Q))$ denotes the O'Neill tensor as above,  $\mathrm{grad}^{\hat{\nabla}}f:=\left(\nabla^{\hat{M}}f\right)^Q$, for every function $f$, and $\xi:=\frac{1}{\rho\sigma}\hat{\xi}$.
\end{elemme}

{\it Proof:}
Since $\beta$ is chosen to be equal to $1$,  the first identity in \eqref{equationwarpedproduct} becomes $\wnabla_{\partial_t}\partial_t=0$ and the other three identities imply that for all $X,Y\in\Gamma(\pi_2^*TM)$ the following relations hold: 
\be
\wnabla_{\partial_t}X&=&\frac{\partial X}{\partial t}+\frac{1}{2}g_t^{-1}\frac{\partial g_t}{\partial t}(X,\cdot),\\
\wnabla_X \partial_t&=&\frac{1}{2}g_t^{-1}\frac{\partial g_t}{\partial t}(X,\cdot),\\
\wnabla_X Y&=&\nabla_X^{M_t} Y-\frac{1}{2}\frac{\partial g_t}{\partial t}(X,Y)\partial_t,
\ee
with $\frac{\partial g_t}{\partial t}=2(\rho\sigma)'(t)(\rho\sigma)(t)\hat{g}_{\hat{\xi}}+2\frac{\partial(\rho k)}{\partial t}(t,x)(\rho(t)k(t,x))\hat{g}_{\hat{\xi}^\perp}$ and hence 
$$g_t^{-1}\frac{\partial g_t}{\partial t}(X,\cdot)=2\frac{(\rho\sigma)'}{\rho\sigma}\hat{g}(X,\hat{\xi})\hat{\xi}\oplus\frac{2}{\rho k}\cdot\frac{\partial(\rho k)}{\partial t}X^\perp=2\frac{(\rho\sigma)'}{\rho\sigma}\wit{g}(X,\xi)\xi\oplus\frac{2}{\rho k}\cdot\frac{\partial(\rho k)}{\partial t}X^\perp,$$
where $\xi=\frac{1}{\rho\sigma}\hat{\xi}$ and $X=\hat{g}(X,\hat{\xi})\hat{\xi}+X^\perp$, with $X^\perp\in\Gamma(\pi_2^*Q)$.
Note that the vector field $\xi$ has unit length with respect to the metric $\wit{g}$. Thus, we obtain the following identities:
\be 
\wnabla_{\partial_t}\partial_t&=&0,\\
\wnabla_{\partial_t}X&=&\frac{\partial X}{\partial t}+\frac{(\rho\sigma)'}{\rho\sigma}\wit{g}(X,\xi)\xi\oplus\frac{1}{\rho k}\cdot\frac{\partial(\rho k)}{\partial t}X^\perp,\\
\wnabla_X \partial_t&=&\frac{(\rho\sigma)'}{\rho\sigma}\wit{g}(X,\xi)\xi\oplus\frac{1}{\rho k}\cdot\frac{\partial(\rho k)}{\partial t}X^\perp,\\
\wnabla_X Y&=&\nabla_X^{M_t} Y-\left(\frac{(\rho\sigma)'}{\rho\sigma}\wit{g}(X,\xi)\wit{g}(Y,\xi)+\frac{1}{\rho k}\cdot\frac{\partial(\rho k)}{\partial t}\wit{g}(X^\perp,Y^\perp)\right)\partial_t.
\ee
Now, we need to compute $\nabla_X^{M_t} Y$ in a more precise way according to the components $X$ and $Y$ in the orthogonal splitting $TM=\mathbb{R}\cdot\xi\oplus Q$.
Recall the Koszul formula, valid for any $X,Y,Z\in\Gamma(TM)$:
\begin{eqnarray}\label{eq:Koszulformula} 
\nonumber g_t(\nabla_X^{M_t} Y,Z)&=&\frac{1}{2}\Big\{X(g_t(Y,Z))+Y(g_t(Z,X))-Z(g_t(X,Y))\\
&&\phantom{\frac{1}{2}\Big\{}+g_t([X,Y],Z)-g_t([Y,Z],X)+g_t([Z,X],Y)\Big\}.
\end{eqnarray}
 First we consider the case when $Y=\xi$.
For $X=\xi$, we have $g_t(\nabla_{\xi}^{M_t}\xi,\xi)=0$ and, for every $Z\in\Gamma(Q)$,
$$g_t(\nabla_{\xi}^{M_t}\xi,Z)=-g_t([\xi,Z],\xi)=-\hat{g}([\hat\xi,Z],\hat\xi)=\hat{g}(Z,\underbrace{\nabla_{\hat{\xi}}^{\hat{M}}\hat{\xi}}_{0})=0,$$
so that $\nabla_{\xi}^{M_t}\xi=0$.
For $X=Z\in\Gamma(Q)$, we have $g_t(\nabla_{Z}^{M_t}\xi,\xi)=0$ and, for every $Z'\in\Gamma(Q)$,
\be 
g_t(\nabla_{Z}^{M_t}\xi,Z')&=&\frac{1}{2}\Big\{Z(g_t(\xi,Z'))+\xi(g_t(Z',Z))-Z'(g_t(Z,\xi))+g_t([Z,\xi],Z')-g_t([\xi,Z'],Z)+g_t([Z',Z],\xi)\Big\}\\
&=&\frac{1}{2}\Big(\xi((\rho k)^2\hat{g}(Z,Z'))-\frac{(\rho k)^2}{\rho\sigma}\hat{g}([\hat{\xi},Z],Z')-\frac{(\rho k)^2}{\rho\sigma}\hat{g}([\hat{\xi},Z'],Z)+(\rho\sigma)\hat{g}([Z',Z],\hat{\xi})\Big)\\
&=&\frac{\xi(k)}{k}g_t(Z,Z')+\frac{(\rho k)^2}{2\rho\sigma}\left(\hat{\xi}(\hat{g}(Z,Z'))-\hat{g}([\hat{\xi},Z],Z')-\hat{g}([\hat{\xi},Z'],Z)\right)+\rho\sigma\hat{g}(\nabla_Z^{\hat{M}}\hat{\xi},Z')\\
&=&\frac{\xi(k)}{k}g_t(Z,Z')+\frac{\sigma}{\rho k^2}g_t(\hat{h}Z,Z'),
\ee
so that $\nabla_{Z}^{M_t}\xi=\frac{\xi(k)}{k}Z+\frac{\sigma}{\rho k^2}\hat{h}Z$.
In the last equality, we used the fact that $\hat\xi$ is a Killing vector field with respect to the metric $\hat{g}$. Let us now choose $Y=Z'\in\Gamma(Q)$ and compute as follows:
\be 
\nabla_{\xi}^{M_t}Z'&=&\nabla_{Z'}^{M_t}\xi-[Z',\xi]\\
&=&\frac{\xi(k)}{k}Z'+\frac{\sigma}{\rho k^2}\hat{h}Z'+\frac{1}{\rho\sigma}[\hat{\xi},Z']\\
&=&\frac{\xi(k)}{k}Z'+\frac{\sigma}{\rho k^2}\hat{h}Z'+\frac{1}{\rho\sigma}[\hat{\xi},Z']^Q+\frac{1}{\rho\sigma}\underbrace{\hat{g}([\hat{\xi},Z'],\hat{\xi})}_{0}\hat{\xi}\\
&=&\frac{\xi(k)}{k}Z'+\frac{\sigma}{\rho k^2}\hat{h}Z'+\hat{\nabla}_{\xi}Z'.
\ee
On the other hand, for every $Z\in\Gamma(Q)$, we have 
$$g_t(\nabla_Z^{M_t}Z',\xi)=-g_t(\nabla_Z^{M_t}\xi,Z')=-\frac{\xi(k)}{k}g_t(Z,Z')-\frac{\sigma}{\rho k^2}g_t(\hat{h}Z,Z')$$ 
and, for any $Z''\in\Gamma(Q)$ we compute:
\be 
g_t(\nabla_Z^{M_t}Z',Z'')&=&\frac{1}{2}\Big\{Z(g_t(Z',Z''))+Z'(g_t(Z'',Z))-Z''(g_t(Z,Z'))\\
&&\phantom{\frac{1}{2}\Big\{}+g_t([Z,Z'],Z'')-g_t([Z',Z''],Z)+g_t([Z'',Z],Z')\Big\}\\
&=&\frac{1}{2}\Big\{Z((\rho k)^2\hat{g}(Z',Z''))+Z'((\rho k)^2\hat{g}(Z'',Z))-Z''((\rho k)^2\hat{g}(Z,Z'))\\
&&\phantom{\frac{1}{2}\Big\{}+(\rho k)^2\hat{g}([Z,Z'],Z'')-(\rho k)^2\hat{g}([Z',Z''],Z)+(\rho k)^2\hat{g}([Z'',Z],Z')\Big\}\\
&=&\rho^2Z(k)k\hat{g}(Z',Z'')+\rho^2Z'(k)k\hat{g}(Z'',Z)-\rho^2Z''(k)k\hat{g}(Z,Z')+(\rho k)^2\hat{g}(\nabla_Z^{\hat{M}}Z',Z'')\\
&=&\frac{Z(k)}{k}g_t(Z',Z'')+\frac{Z'(k)}{k}g_t(Z'',Z)-\frac{Z''(k)}{k}g_t(Z,Z')+g_t(\hat{\nabla}_ZZ',Z'').
\ee
This finishes the proof of the lemma. 
\findemo

Let us now recall some standard definitions. On a Riemannian flow $(M,\hat{g},\hat{\xi})$, a function $f$ is said to be {\it basic} if $\hat{\xi}(f)=0$, \emph{i.e.} the function $f$ depends only on the transversal variables. A {\it transversal} K\"ahler structure $J$ on a Riemannian flow is defined as an almost-Hermitian structure $J:\Gamma(Q)\to \Gamma(Q)$, which is parallel with respect to the transversal Levi-Civita connection. The following lemma will be useful when considering basic transversal conformal changes of the metric on a Riemannian flow.  Indeed, we will show that for any conformal change of the transverse metric by a basic function, the flow will be still Riemannian. 

\begin{elemme}\label{l:transvconfchange}
Let $(M,\hat{g},\hat{\xi})$ be a connected minimal Riemannian flow and let $f\in C^\infty(M,\R)$ be a basic function. 
Then $(M,g:=\hat{g}_{\hat{\xi}}\oplus e^{2f}\cdot\hat{g}_{\hat{\xi}^\perp},\xi:=\hat{\xi})$ is a minimal Riemannian flow, whose O'Neill tensor is given by $h=e^{-2f}\hat{h}$ and whose Levi-Civita connection satisfies for all $X\in\Gamma(TM)$ and $Z\in\Gamma(Q)$:
$$\nabla_XZ=\hat{\nabla}_XZ+X^Q(f)Z+Z(f)X^Q-\hat{g}(X,Z)\hat{\nabla}f,$$
where $\hat{\nabla}f:=\mathrm{grad}^{\hat{\nabla}}f$ is the pointwise projection of the $\hat{g}$-gradient of $f$ onto $Q$. Moreover, if $J$ is a transversal K\"ahler structure on $(M,\hat{g},\hat{\xi})$, then $J$ remains a transversal K\"ahler structure on $(M,g,\xi)$ if and only if either $\mathrm{rk}(Q)=2$ or $f$ is constant.
\end{elemme}

{\it Proof:} First, we make use of the Koszul formula (\ref{eq:Koszulformula}) to show that for any $Z\in\Gamma(Q)$ we have:
$$g(\nabla_\xi^M\xi,Z)=-g([\xi,Z],\xi)=-\hat{g}([\hat{\xi},Z],\hat{\xi})=\hat{g}(\nabla_{\hat{\xi}}^{\hat{M}}\hat{\xi},Z)=0.$$
Moreover, the Lie derivative of the transverse conformal metric in the direction vector field $\hat{\xi}$ is equal to $$\mathcal{L}_\xi(e^{2f}\cdot\hat{g}_{\hat{\xi}^\perp})=e^{2f}\xi(f)\hat{g}_{\hat{\xi}^\perp}+e^{2f}\mathcal{L}_\xi(\hat{g}_{\hat{\xi}^\perp})=0,$$
since $f$ is assumed to be a basic function. In particular, this shows that $\xi$ is a unit Killing vector field with respect to the metric $g$ and therefore $(M,g,\xi)$ is a minimal Riemannian flow. The relation between $\nabla$ and $\hat{\nabla}$ is proven as in the usual case by the uniqueness of a compatible transversal torsion-free connection. To compare the corresponding O'Neill tensors, we just compute for any $Y,Z\in \Gamma(Q)$ as follows:
$$g(hY,Z)=-\frac{1}{2}g_\xi([Y,Z],\xi)=-\frac{1}{2}\hat{g}_{\hat{\xi}}([Y,Z],\hat{\xi})=\hat{g}(\hat{h}Y,Z)=e^{-2f}g(\hat{h}Y,Z).$$
Let $J$ be a transversal K\"ahler structure on $(M,\hat{g},\hat{\xi})$. Then, $J$ remains an almost-Hermitian structure on $Q$ and $\nabla J=0$ if and only if, for all $Z,Z'\in\Gamma(Q)$, the following equality holds:
$$Z(f)JZ'+JZ'(f)Z-g(Z,JZ')\hat{\nabla}f=Z(f)JZ'+Z'(f)JZ-g(Z,Z')J\hat{\nabla}f,$$
which is equivalent to
$$Z'(f)JZ-JZ'(f)Z+g(Z,JZ')\hat{\nabla}f-g(Z,Z')J\hat{\nabla}f=0.$$
In case $\mathrm{rk}(Q)=2$, this identity is trivially satisfied for all $Z,Z'\in\Gamma(Q)$, whereas for $\mathrm{rk}(Q)>2$ it is satisfied  if and only if $\hat{\nabla}f=0$, that is, $f$ is constant.
This concludes the proof of Lemma~\ref{l:transvconfchange}. \findemo

\begin{erem}\label{rem:leafwisechange}
{\rm Note that if we rescale the vector field $\hat{\xi}$ by some positive number, that is, we consider the metric $g:=\alpha^2\hat{g}_{\hat{\xi}}\oplus \hat{g}_{\hat{\xi}^\perp}$ for some positive $\alpha$, then $\left(M,g,\xi:=\frac{1}{\alpha}\hat{\xi}\right)$ is still a minimal Riemannian flow with O'Neill tensor $h=\alpha\hat{h}$. }
\end{erem}
\color{black}

 Let $(M,\hat{g},\hat{\xi})$ be a minimal Riemannian flow and assume the existence of a transversal K\"ahler structure $J$ on $Q=\hat{\xi}^\perp$. We consider the  almost-Hermitian structure $\wit{J}$ on $\wit{M}=I\times M$ defined by setting 
\begin{equation}\label{newcomplexstructure}
\wit{J}(\partial_t):=-\frac{1}{\rho\sigma}\hat{\xi},\,\, \wit{J}(\frac{1}{\rho\sigma}\hat{\xi}):=\partial_t\,\, \text{and}\,\, \wit{J}_{|_{\{\hat{\xi},\partial_t\}^\perp}}:=J,
\end{equation}
where $\rho,\sigma,k$ are the coefficients of the metric in \eqref{formofthemetric}. Similarly to \cite[Lemma 3.4]{GinSemm11}, we want to characterize in the next lemma those functions $\rho,\sigma,k$ for which $(\wit{M}^{2n},\wit{g},\wit{J})$ is \emph{K\"ahler}, that is, for which $\wnabla\wit{J}=0$ holds on $\wit{M}$. Note in particular that $M$ has odd dimension equal to $2n-1$.

\begin{elemme}\label{l:charactKaehlerrhosigmak}
Let $(\wit{M},\wit{g}):=\left(I\times M,dt^2\oplus\rho(t)^2\left(\sigma(t)^2\hat{g}_{\hat{\xi}}\oplus k(t,x)^2\hat{g}_{\hat{\xi}^\perp}\right)\right)$ be a doubly-warped product, where $(M,\hat{g},\hat{\xi})$ is a minimal Riemannian flow and $\rho,\sigma\colon I\to\R_+^\times$, as well as $k\colon I\times M\to\R_+^\times$, are smooth positive functions. Let $\xi:=\frac{1}{\rho\sigma}\hat \xi$ and $\hat h :=\nabla^{\hat M}\hat{\xi}$. Furthermore, we assume $(M,\hat{g},\hat{\xi})$ carries a transversal K\"ahler structure $J$ and we define the almost-Hermitian structure $\wit{J}$ on $\wit{M}$ via \eqref{newcomplexstructure}. Then the following statements hold true:
\begin{enumerate}[i)]
\item The structure $(\wit{M}^{2n},\wit{g},\wit{J})$ is K\"ahler if and only if $\hat{\xi}(k)=0$, $\hat{h}=-\frac{k}{\sigma}\cdot\frac{\partial(\rho k)}{\partial t}J$ and, if $n>2$, $\mathrm{grad}^{\hat{\nabla}}(k)=0$ (and thus the function $k$ only depends on $t$, if $M$ is connected).
In this case, there exists a basic function $C$ on $M$, which is constant if $n>2$, such that $\frac{\partial(\rho k)^2}{\partial t}=2\rho\sigma C$.
\item If $\hat{h}$ vanishes at a point (for $n>2$) or vanishes identically (for $n=2$), then $(\wit{M}^{2n},\wit{g},\wit{J})$ is K\"ahler if and only if it is locally isometric to a K\"ahler product $(\R_t\times\R_s\times \Sigma, dt^2\oplus \rho^2(t)ds^2\oplus g_{\Sigma})$, for some K\"ahler manifold $(\Sigma, g_{\Sigma})$ and some positive function $\rho$ on $I$ (that plays the role of $\rho\sigma$).
\item If $(\wit{M}^{2n},\wit{g},\wit{J})$ is K\"ahler, $\hat{h}\neq0$ on $M$ and $k$ only depends on $t$, then, up to rescaling $\hat{\xi}$, turning $t$ into $-t$ or setting $\ovl{\rho} :=\rho k$, as well as $\ovl{\sigma}:=\frac{\sigma}{k}$, we have $\rho'=\sigma$ on $I$ and $k=1$, hence $\hat{h}=-J$ on $M$. In particular, $(M,\hat{g},\hat{\xi})$ is Sasaki.
\item\label{point:kproduct}  If $(\wit{M}^{2n},\wit{g},\wit{J})$ is K\"ahler, $\hat{h}\neq0$ on $M$ and $k$ is of the form $k(t,x)=k_1(t)k_2(x)$ (hence $n=2$), then, up to turning $t$ into $-t$, rescaling $\hat{\xi}$, setting $\ovl{\rho} :=\rho k_1$, as well as $\ovl{\sigma}:=\frac{\sigma}{k_1}$, we may assume that $k_1=1$ and there exists a basic positive function $C$ on $M$ such that $\rho'=\sigma$ on $I$ and $k_2=\sqrt{C}$, hence $\hat{h}=-C\cdot J$ on $M$. 
In particular, $\left(M^3,\bg:=\hat{g}_{\hat{\xi}}\oplus C\cdot\hat{g}_{\hat{\xi}^\perp},\hat{\xi}\right)$ is Sasaki. In this case, we call $(M,\hat{g},\hat{\xi})$ basic conformally Sasaki.
\end{enumerate}
\end{elemme}

{\it Proof:} \begin{enumerate}[$i)$]
\item  We first compute $\wnabla\wit{J}$ using Lemma \ref{l:identitiesrhosigmak}.
We keep the same notations as in the proof of Lemma \ref{l:identitiesrhosigmak}.
Let $Z,Z'\in\Gamma(\pi_2^*Q)$ be arbitrary sections.
First,
$$\wnabla_{\partial_t}(\wit{J}\partial_t)-\wit{J}(\wnabla_{\partial_t}\partial_t)=-\wnabla_{\partial_t}\xi+0=0$$
as well as $\wnabla_{\partial_t}(\wit{J}\xi)-\wit{J}(\wnabla_{\partial_t}\xi)=0$.
Moreover, 
\be 
\wnabla_{\partial_t}(\wit{J}Z)-\wit{J}(\wnabla_{\partial_t}Z)&=&\frac{\partial JZ}{\partial t}+\frac{1}{\rho k}\cdot\frac{\partial(\rho k)}{\partial t}JZ-\wit{J}\left(\frac{\partial Z}{\partial t}+\frac{1}{\rho k}\cdot\frac{\partial(\rho k)}{\partial t}Z\right)\\
&=&\left[\frac{\partial}{\partial t},J\right]Z\\
&=&0,
\ee
showing that $\wnabla_{\partial_t}\wit{J}=0$. Now, differentiating in the direction of $\xi$, we first obtain
$$\wnabla_\xi(\wit{J}\partial_t)-\wit{J}(\wnabla_\xi\partial_t)=-\wnabla_\xi\xi-\frac{(\rho\sigma)'}{\rho\sigma}\wit{J}\xi=0$$
as well as $\wnabla_\xi(\wit{J}\xi)-\wit{J}(\wnabla_\xi\xi)=0$.
Furthermore, using $\hat{\nabla}J=0$,
\be 
\wnabla_\xi(\wit{J}Z)-\wit{J}(\wnabla_\xi Z)&=&\hat{\nabla}_{\xi}(JZ)+\frac{\xi(k)}{k}JZ+\frac{\sigma}{\rho k^2}\hat{h}JZ-\wit{J}\left(\hat{\nabla}_{\xi}Z+\frac{\xi(k)}{k}Z+\frac{\sigma}{\rho k^2}\hat{h}Z\right)\\
&=&\frac{\sigma}{\rho k^2}\left[\hat{h},J\right]Z.
\ee
Thus, $\wnabla_\xi\wit{J}=0$ if and only if $\left[\hat{h},J\right]=0$ on $M$.
It remains to look at differentiation in transversal directions:
\be 
\wnabla_Z(\wit{J}\partial_t)-\wit{J}(\wnabla_Z\partial_t)&=&-\wnabla_Z\xi-\wit{J}\left(\frac{1}{\rho k}\cdot\frac{\partial(\rho k)}{\partial t}Z\right)\\
&=&-\frac{\xi(k)}{k}Z-\frac{\sigma}{\rho k^2}\hat{h}Z-\frac{1}{\rho k}\cdot\frac{\partial(\rho k)}{\partial t}JZ,
\ee
in particular $(\wnabla_Z\wit{J})\partial_t=0$ for all $Z\in\Gamma(\pi_2^*Q)$ if and only if $\xi(k)=0$ and $\hat{h}=-\frac{k}{\sigma}\cdot\frac{\partial\rho k}{\partial t}J$; note that the latter condition implies $\left[\hat{h},J\right]=0$.
Similarly, 
\be 
\wnabla_Z(\wit{J}\xi)-\wit{J}(\wnabla_Z\xi)&=&\frac{1}{\rho k}\cdot\frac{\partial(\rho k)}{\partial t}Z-\wit{J}\left(\frac{\xi(k)}{k}Z+\frac{\sigma}{\rho k^2}\hat{h}Z\right)\\
&=&-\frac{\xi(k)}{k}JZ+\frac{1}{\rho k}\cdot\frac{\partial(\rho k)}{\partial t}Z-\frac{\sigma}{\rho k^2}J\hat{h}Z,
\ee
so that $(\wnabla_Z\wit{J})\xi=0$ for all $Z\in\Gamma(\pi_2^*Q)$ if and only if $\xi(k)=0$ and $\hat{h}=-\frac{k}{\sigma}\cdot\frac{\partial(\rho k)}{\partial t}J$, which is precisely what we had before.
Last but no the least, assuming the latter conditions are fulfiled, we obtain
\be 
\wnabla_Z(\wit{J}Z')-\wit{J}(\wnabla_ZZ')&=&\hat{\nabla}_Z(JZ')+Z(\ln(k))JZ'+JZ'(\ln(k))Z-\wit{g}(Z,JZ')\mathrm{grad}^{\hat{\nabla}}(\ln(k))\\
&&-\xi(\ln(k))\wit{g}(Z,JZ')\xi-\frac{\sigma}{\rho k^2}\wit{g}(\hat{h}Z,JZ')\xi-\frac{\partial\ln(\rho k)}{\partial t}\wit{g}(Z,JZ')\partial_t\\
&&-J(\hat{\nabla}_ZZ')-Z(\ln(k))JZ'-Z'(\ln(k))JZ+\wit{g}(Z,Z')J(\mathrm{grad}^{\hat{\nabla}}(\ln(k)))\\
&&+\xi(\ln(k))\wit{g}(Z,Z')\partial_t+\frac{\sigma}{\rho k^2}\wit{g}(\hat{h}Z,Z')\partial_t-\frac{\partial\ln(\rho k)}{\partial t}\wit{g}(Z,Z')\xi\\
&=&JZ'(\ln(k))Z-Z'(\ln(k))JZ+\wit{g}(JZ,Z')\mathrm{grad}^{\hat{\nabla}}(\ln(k))\\&&+\wit{g}(Z,Z')J(\mathrm{grad}^{\hat{\nabla}}(\ln(k))).
\ee
We now want to know when the r.h.s. of the last identity vanishes.
In case $2n-2=2$, i.e. $n=2$, it trivially vanishes pointwise for all $Z,Z'\in Q$.
In case $n>2$, assuming $\mathrm{grad}^{\hat{\nabla}}(\ln(k))$ to be nonzero at a point, we may chose for instance $Z\in Q$ such that $Z,JZ,\mathrm{grad}^{\hat{\nabla}}(\ln(k)),J(\mathrm{grad}^{\hat{\nabla}}(\ln(k)))$ are linearly independent, in which case it can be deduced from the vanishing of all terms in the r.h.s. of the last identity that $\wit{g}(JZ,Z')=0$ for all $Z'\in Q$, which is a contradiction; thus $\mathrm{grad}^{\hat{\nabla}}(\ln(k))$ must vanish identically for $n>2$. To sum up, we have shown that $\wnabla\wit{J}=0$ on $\wit{M}$ if and only if $\xi(k)=0$, $\hat{h}=-\frac{k}{\sigma}\cdot\frac{\partial\rho k}{\partial t}J$ and, if $n>2$, $\mathrm{grad}^{\hat{\nabla}}(k)=0$ -- and hence $k$ only depends on $t$, if $M$ is connected.
\item If $\hat h=0$ at a point if $n>2$ or identically if $n=2$ and if $(\tilde M,\tilde g, \tilde J)$ is K\"ahler, then by $i)$ $\hat h$ must vanish identically and the following identity $\frac{\partial(\rho k)}{\partial t}=0$ holds.
Hence, $\rho k$ is a function depending only on $x$: $\rho(t)k(t,x)=D(x)$.
If $n>2$, as $k$ depends only on $t$ as we have already shown, then $D$ is constant and in this case we have: $\tilde g=dt^2\oplus\rho^2(t)\sigma^2(t)\hat g_{\hat \xi}\oplus D^2 \hat g_{\hat \xi^{\perp}}$.
Rescaling $\hat g_{\hat{\xi}^\perp}$ and replacing $\rho\sigma$ by $\rho$, we obtain the desired product form.
If $n=2$, then $D$ may be nonconstant, but $\hat\xi(D)=0$ and in this case $(\Sigma, g_{\Sigma})$ is a surface (hence any Hermitian metric on $\Sigma$ is already K\"ahler).
\item If $k$ only depends on $t$ and $(\tilde{M},\tilde{g},\tilde{J})$ is K\"ahler then we may assume, up to replacing $\rho$ by $\rho k$ and $\sigma$ by $\frac{\sigma}{k}$, that $k=1$.
Since neither $\hat{h}$ nor $J$ depend on $t$, there exists a constant $C$ such that $\frac{\rho'}{\sigma}=C$.
If $C=0$, then $\hat h=0$ and we are back in case $ii)$.
If $C\neq0$, then up to rescaling $\hat{\xi}$, we may assume that $C=\pm1$ and, up to turning $t$ into $-t$, that $C=1$.
Then $\rho'=\sigma$ on $I$ and $\hat{h}=-J$ on $M$.
\item In case $k$ is not necessarily constant in $x$ (and thus $n=2$), there exists a function $C$ on $M$, which must be basic since both $\hat{h}$ and $J$ are, such that $\frac{k}{\sigma}\cdot\frac{\partial\rho k}{\partial t}=C$, that is, $\frac{\partial(\rho k)^2}{\partial t}=2\rho\sigma\cdot C$ on $\wit{M}$.
This is equivalent to $(\rho(t)k(t,x))^2-(\rho(0)k(0,x))^2=2C(x)\cdot\int_0^t\rho(s)\sigma(s)\,ds$ for all $(t,x)\in I\times M$.
This shows that $k^2$ is the sum of two functions that are products of a function of $t$ with a function of $x$, still $k$ must not be itself in product form.
In case $k$ is of the form $k(t,x)=k_1(t)k_2(x)$, then we may assume as above that $k_1=1$ (up to changing $\rho$ and $\sigma$ by multiplying $\rho$ by $k_1$ and $\sigma$ by $\frac{1}{k_1}$).
The identity $\frac{\partial(\rho k_2)^2}{\partial t}=2\rho\sigma\cdot C$ yields $\rho' k_2^2=\sigma\cdot C$. If $\rho'$ vanishes at one point, then $C$ must vanish identically and then $\hat{h}=0$, which is again the case $ii)$.
Otherwise, up to turning $t$ into $-t$, we may assume $\rho'>0$ on $I$, from which $k_2(x)^2=\frac{\sigma}{\rho'}(t)\cdot C(x)$ follows, in particular $\frac{\rho'}{\sigma}$ is constant.
Up to rescaling $\hat{\xi}$, we may assume that $\rho'=\sigma$ on $I$, from which $k_2(x)^2=C(x)$ follows for all $x\in M$ (showing on the way that $C$ must be positive) and therefore $k_2=\sqrt{C}$. The last claim follows from Lemma \ref{l:transvconfchange}.
This concludes the proof of Lemma~\ref{l:charactKaehlerrhosigmak}.
\end{enumerate}

 \begin{erem}
 	{\rm The product form assumed for the function $k$ in Lemma~\ref{l:charactKaehlerrhosigmak}, $\ref{point:kproduct})$, is fulfiled in the case where we apply this result, see Theorem~\ref{t:cardwp2eigenv}.}
\end{erem}

We end this section by characterizing the completeness of doubly-warped products. We will consider the case when $\mathrm{grad}^{\hat{\nabla}}k=0$ (or $k=1$, up to rescaling the metric) in the expression \eqref{formofthemetric} of the metric, since these cases naturally arise in our classification results (see Theorems \ref{t:caracdwpEinstein} and \ref{t:cardwp2eigenv}). 

\color{black}


\begin{elemme}\label{l:completedwp} 
Let $(\wit{M},\wit{g}):=\left(I\times M,dt^2\oplus\rho(t)^2(\sigma(t)^2\hat{g}_{\hat{\xi}}\oplus\hat{g}_{\hat{\xi}^\perp})\right)$ be a doubly-warped product. Then $(\wit{M},\wit{g})$ is complete if and only if $(M,\hat{g})$ is complete and $I=\R$. 
\end{elemme} 

{\it Proof:} The proof follows that of the analogous result for warped products \cite[Lemma 7.2]{BishopONeill69}. Assume $(\wit{M},\wit{g})$ to be complete. Then $(M,\hat{g})$ must complete because it is a closed subset of $\wit{M}$ and the metrics $g_t:=\rho(t)^2(\sigma(t)^2\hat{g}_{\hat{\xi}}\oplus\hat{g}_{\hat{\xi}^\perp})$ and $\hat{g}$ are equivalent for any fixed $t$. Moreover, because the integral curves of $\partial_t$ are geodesics, according to \eqref{equationwarpedproduct}, then $I=\R$ must hold. Conversely, assume $(M,\hat{g})$ to be complete and $I=\R$. Let $\left((t_n,x_n)\right)_n$ be any Cauchy sequence in $(\wit{M},\wit{g})$, then because the distance between projected $\R$-components is anyway smaller than the distance associated to $\wit{g}$, the sequence $(t_n)_n$ must be a Cauchy sequence in $\R$ and therefore must converge to some $T\in\R$. But since furthermore the sequence $(t_n)_n$ must be bounded, so must be the coefficients $\rho(t_n),\sigma(t_n)$ of the metric $g_{t_n}$ independently of $n$, therefore all $(g_{t_n})_n$ are uniformly equivalent to $\hat{g}$. As a consequence, the sequence $(x_n)_n$ must be a Cauchy sequence on $(M,\hat{g})$ and therefore must converge to some $x\in M$. Because again $(g_{t_n})_n$ are uniformly equivalent to $\hat{g}$, the sequence $\left((t_n,x_n)\right)_n$ must converge to $(T,x)\in\wit{M}$, which concludes the proof. \findemo

\section{Proof of the main result}\label{s:Obatadwp}



In this section, we give the proof of Theorem \ref{t:cardwp2eigenv}. For this, we let $\nu:=\frac{\wnabla u}{|\wnabla u|}\in\Gamma(T\wit{M})$ and $\xi:=-\wit{J}\nu$. By assumption,
$$\wnabla^2u=\lambda\cdot\left(\nu^\flat\!\otimes\!\nu\oplus\xi^\flat\!\otimes\!\xi\right)+\mu\cdot\mathrm{Id}_{\{\xi,\nu\}^\perp}$$
on $\wit{M}$ for smooth real-valued functions $\lambda,\mu$ on $\wit{M}$. Fix any value $u_0$ of $u$ and let $M:=u^{-1}(\{u_0\})$, which is a real hypersurface of $\wit{M}$ with induced metric $g:=\wit{g}_{|_M}$.
Since by assumption $\wnabla u$ is a pointwise eigenvector for $\wnabla^2 u$, the vector field $\nu$ is actually geodesic on $(\wit{M},\wit{g})$ and Proposition \ref{p:warpedprodu} shows that the map $F\colon \R\times M\to \wit M$, given by the flow of $\nu$, is a diffeomorphism pulling $\wit{g}$ back onto $dt^2\oplus g_t$, where $g_t:=(F_t)^*g_{|_{TM\times TM}}$ is a smooth one-parameter-family of Riemannian metrics on $M$ that we next determine more precisely. In the proof of Proposition \ref{p:warpedprodu}, we show that $f:=u\circ F$ only depends on $t$. Let $t_0$ be such that $f(t_0)=(u\circ F)(t_0)=u_0$. Moreover, we have $f'(t)=|\wnabla u|_{F(t,x)}$ and, since $\nu$ is geodesic, $f''(t)=\lambda\circ F(t,x)$. Therefore, as $[\xi,\nu]=\frac{\lambda}{|\wnabla u|}\xi$, which can be shown by a straightforward computation, we have 
\be 
\frac{\partial}{\partial s}((F_s)_*\xi)_{F_{t_0}(x)}{|_{s=t}}&=&\left((F_t)_*\left(\frac{\lambda}{|\wnabla u|}\xi\right)\right)_{F_{t_0}(x)}\\
&=&\frac{\lambda}{|\wnabla u|}\circ F_{-t}(F_{t_0}(x))\cdot((F_t)_*\xi)_{F_{t_0}(x)}\\
&=&\frac{\lambda}{|\wnabla u|}\circ F_{-t+t_0}(x)\cdot((F_t)_*\xi)_{F_{t_0}(x)}\\
&=&\frac{f''(-t+t_0)}{f'(-t+t_0)}\cdot((F_t)_*\xi)_{F_{t_0}(x)}.
\ee
Integrating, we obtain
\begin{equation}\label{eq:differentialchi}
\begin{split}
((F_t)_*\xi)_{F_{t_0}(x)}=&\exp\left(\int_0^t\frac{f''(-s+t_0)}{f'(-s+t_0)}ds\right)\cdot\xi(F_{t_0}(x))\\
=&\exp\left(\int_{-t_0}^{t-t_0}\frac{f''(-s')}{f'(-s')}ds'\right)\cdot\xi(F_{t_0}(x))\textrm{ where }s':=s-t_0\\
=&\frac{f'(t_0)}{f'(t_0-t)}\xi(F_{t_0}(x)).
\end{split}
\end{equation}
In particular, this leads to $(F_t^*\wit{g})(\xi,\xi)=\frac{f'(t)^2}{f'(0)^2}$. On the other hand, we have for all $X,Y\in TM$  
\be 
(\mathcal{L}_\nu \wit{g})(X,Y)&=&2\wit{g}(\wnabla_X\nu,Y)=\frac{2}{|\wnabla u|}\wnabla^2u(X,Y)\\
&=&\frac{2}{|\wnabla u|}\left((\lambda-\mu)\wit{g}(X,\xi)\wit{g}(Y,\xi)+\mu \wit{g}(X,Y)\right).
\ee
Hence as in Proposition \ref{p:warpedprodu}, we get for all $X,Y\in\{\xi,\nu\}^\perp$ that
$$\frac{\partial}{\partial s}F_s^*\wit{g}(X,Y)_{|_{s=t}}=\frac{2}{f'(t)}\mu\circ F(t,x)\cdot(F_t^*\wit{g})(X,Y),$$
in particular after integrating, we find
$$(F_t^*\wit{g})(X,Y)=\exp\left(2\int_0^t\frac{\mu\circ F(s,x)}{f'(s)}\,ds\right)\cdot \wit{g}(X,Y).$$
Therefore,
\begin{equation}\label{eq:Fstarg}
F_t^*\wit{g}=dt^2\oplus\frac{(f')^2(t)}{(f')^2(0)}g_{\xi}\oplus\exp\left(2\int_0^t\frac{\mu\circ F(s,x)}{f'(s)}\,ds\right) g_{\xi^\perp},
\end{equation}
where we recall that $g$ is the induced metric $g=\wit{g}_{|_M}$ and $g=g_{\xi}\oplus g_{\xi^\perp}$. We check now that $(M,g,\xi)$ is a minimal Riemannian flow. Firstly, since for all $X\in T\wit{M}$
\begin{equation} \label{eq:derivativenablaxu}
\wnabla_X\nu=\frac{1}{|\wnabla u|}\left(\wnabla_X^2u-g(\wnabla_X^2u,\nu)\nu\right)=\frac{1}{|\wnabla u|}\left(\wnabla_X^2u-\lambda g(X,\nu)\nu\right),
\end{equation}
we have 
$$\wnabla_\xi\xi=\wnabla_{\wit J\nu}\wit J\nu=\wit J\left(\wnabla_{\wit J\nu}\nu\right)=-\frac{1}{|\wnabla u|}\wnabla_{\nu}^2u=-\frac{\lambda}{|\wnabla u|}\nu$$
and hence $\nabla_\xi^M\xi=0$. 
Secondly, for all $X\in TM\cap\xi^\perp$, we have
\begin{eqnarray} \label{ONeilltens}
\nabla_X^M\xi&=&\wnabla_X\xi+\underbrace{g(\wnabla_X\nu,\xi)}_{0}\nu\nonumber\\
&=&-\wit J(\wnabla_X\nu)\nonumber\\
&=&-\frac{\mu}{|\wnabla u|}\wit JX.
\end{eqnarray}
Thus $\nabla^M\xi$ is skew symmetric and vanishes on $\xi$, therefore $(M,g,\xi)$ is a minimal Riemannian flow and its O'Neill tensor is given by $h= -\frac{\mu}{|\wnabla u|} \wit J.$ In the following, we denote by $\nabla$ the transversal Levi-Civita connection of the flow. Hence both connections are related by $\nabla_\xi Z=\nabla_\xi^MZ-\nabla_Z^M\xi$ and $\nabla_{Z'}Z=\nabla_{Z'}^MZ+g(\nabla_{Z'}^M\xi,Z)\xi$ for all sections $Z,Z'$ of $Q:=\{\xi,\nu\}^\perp\to M$.
As for $J:=\wit{J}$ on $Q$, we have, for every $Z\in\Gamma(Q)$, 
\be 
\nabla_\xi(JZ)&=&\wnabla_\xi(\wit{J}Z)+\underbrace{g(\wnabla_\xi\nu,JZ)}_{0}\nu-\nabla_{JZ}^M\xi\\
&=&\wit{J}(\wnabla_\xi Z)-\nabla_{JZ}^M\xi\\
&=&\wit{J}(\wnabla_\xi Z-\nabla_Z^M\xi)\qquad\textrm{since }[\nabla^M\xi,J]=0\\
&=&\wit{J}(\nabla_\xi^MZ-\nabla_Z^M\xi)\qquad\textrm{since }g(\wnabla_\xi Z,\nu)=0\\
&=&J(\nabla_\xi Z)
\ee
and, for all $Z'\in Q$,
\be 
\nabla_{Z'}(JZ)&=&\wnabla_{Z'}(\wit{J}Z)+g(\wnabla_{Z'}\nu,JZ)\nu+g(\nabla_{Z'}^M\xi,JZ)\xi\\
&=&\wit{J}(\wnabla_{Z'}Z)+\frac{\mu}{|\wnabla u|}g(Z',JZ)\nu-\frac{\mu}{|\wnabla u|}g(JZ',JZ)\xi\\
&=&\wit{J}\Big(\wnabla_{Z'}Z+\frac{\mu}{|\wnabla u|}g(Z',Z)\nu-\frac{\mu}{|\wnabla u|}g(JZ',Z)\xi\Big)\\
&=&J(\nabla_{Z'}Z).
\ee
Therefore, $\nabla J=0$ and hence $J$ defines a transversal K\"ahler structure on $(M,g,\xi)$.\\
In the following, we show that the pull-back of the almost complex structure $F^*\wit{J}$ on $\R\times M$ maps $\xi$ onto $\partial_t$, $\partial_t$ onto $-\xi$ and coincides with $J$ on $Q$.
For all $(t,x)\in\R\times M$ and $X\in\R\oplus T_xM$, we have
$$(F^*\wit{J})_{(t,x)}(X)=\left(d_{(t,x)}F\right)^{-1}\circ\wit{J}_{F(t,x)}\circ\left(d_{(t,x)}F\right)(X).$$
For $X=\partial_t$, we have 
\be 
(F^*\wit{J})_{(t,x)}(\partial_t)&=&\left(d_{(t,x)}F\right)^{-1}\circ\wit{J}_{F(t,x)}\circ\left(d_{(t,x)}F\right)(\partial_t)\\
&=&\left(d_{(t,x)}F\right)^{-1}\circ\wit{J}_{F(t,x)}(\partial_t)\\
&=&-\left(d_{(t,x)}F\right)^{-1}(\xi_{F(t,x)})\\
&=&-(F_{-t}{}_*\xi)_x\\
&\bui{=}{\eqref{eq:differentialchi}}&-\frac{f'(0)}{f'(t)}\xi_x=-\xi_{(t,x)}
\ee
as well as 
\be
(F^*\wit{J})_{(t,x)}(\xi_{(t,x)})&=&\left(d_{(t,x)}F\right)^{-1}\circ\wit{J}_{F(t,x)}\circ\left(d_{(t,x)}F\right)(\xi_{(t,x)})\\
&=&\frac{f'(0)}{f'(t)}\left(d_{(t,x)}F\right)^{-1}\circ\wit{J}_{F(t,x)}\left(\left(d_{(t,x)}F\right)(\xi_{x})\right)\\
&=&\frac{f'(0)}{f'(t)}\left(d_{(t,x)}F\right)^{-1}\circ\wit{J}_{F(t,x)}\left((F_t{}_*\xi)_{F_t(x)}\right)\\
&\bui{=}{\eqref{eq:differentialchi}}&\frac{f'(0)}{f'(t)}\left(d_{(t,x)}F\right)^{-1}\circ\wit{J}_{F(t,x)}\left(\frac{f'(t)}{f'(0)}\cdot\xi_{F(t,x)}\right)\\
&=&\left(d_{(t,x)}F\right)^{-1}(\partial_t)=\partial_t.
\ee
To show that $(F^*\wit{J})_{|_Q}=J$, we compute the Lie derivative of $\wit{J}$ in the $\nu$-direction: for every $X\in T\wit{M}$,
\be 
\left(\mathcal{L}_\nu\wit{J}\right)X&=&[\nu,\wit{J}X]-\wit{J}[\nu,X]\\
&=&\wit{J}(\wnabla_\nu X)-\wnabla_{\wit{J}X}\nu-\wit{J}[\nu,X]\\
&=&\wit{J}\wnabla_X\nu-\wnabla_{\wit{J}X}\nu\\
&\bui{=}{\eqref{eq:derivativenablaxu}}&\frac{1}{|\wnabla u|}\left(\wit{J}\wnabla_X^2u-\lambda g(X,\nu)\wit{J}\nu\right)-\frac{1}{|\wnabla u|}\left(\wnabla_{\wit{J}X}^2u-\lambda g(\wit{J}X,\nu)\nu\right)\\
&=&\frac{\lambda}{|\wnabla u|}\cdot\left(g(X,\nu)\xi+g(X,\xi)\nu\right),
\ee
therefore $\mathcal{L}_\nu\wit{J}=\frac{\lambda}{|\wnabla u|}\cdot\left(\nu^\flat\!\otimes\!\xi+\xi^\flat\!\otimes\!\nu\right)$.
Now, for any $(t,x)\in\R\times M$ and $Z\in T_xM\cap\hat{\xi}^\perp$,
\be
\frac{\partial}{\partial s}(F_s^*\wit{J})(Z_x)_{|_{s=t}}&=&\frac{\partial}{\partial s}(F_s^*\wit{J})_{|_{s=t}}(Z_x)\\
&=&\left(\mathcal{L}_\nu\wit{J}\right)(d_xF_t(Z_x))\\
&=&0
\ee
because $F_t$ preserves $Q$.
We deduce that $(F_t^*\wit{J})(Z_x)=\wit{J}_x(Z_x)=JZ_x$, therefore $(F^*\wit{J})_{|_Q}=J$, as claimed.\\
 Summing up, we have shown  that, on the product manifold $I\times M$, the metric $F_t^*\wit{g}$ is determined by \eqref{eq:Fstarg} and the complex structure $F_t^*\wit{J}$ has the form as in \eqref{newcomplexstructure}. Moreover the manifold $(I\times M, F_t^*\wit{g}, F^*\wit{J})$ is K\"ahler and the triple $(M, g, \xi)$ is a minimal Riemannian flow equipped with a transversal K\"ahler structure $J=\wit J$. Recall that $g=g_\xi\oplus g_{\xi^\perp}$ is the induced metric $\wit{g}_{|_M}$. In the following, we will apply Lemma \ref{l:charactKaehlerrhosigmak} in order to obtain the classification result.\\

We begin with the case when $n>2$. We write $F_t^*\wit{g}$ as a doubly-warped product in the following way 
$$F_t^*\wit{g}=dt^2\oplus\rho^2(t)\left(\sigma^2(t) g_{\xi}\oplus k^2(t,x) g_{\xi^\perp}\right),$$ 
where $\rho,\sigma,k$ are positive smooth functions that satisfy the system:
\begin{equation*}
\left\{\begin{array}{ll}
\rho^2(t)\sigma^2(t)&=\frac{(f')^2(t)}{(f')^2(0)}\\
\rho^2(t)k^2(t,x)&=\exp\left(2\displaystyle\int_0^t\frac{\mu\circ F(s,x)}{f'(s)}\,ds\right)\end{array}\right.
\end{equation*}
for all $(t,x)\in\R\times M$. Therefore, as the hypothesis of Lemma \ref{l:charactKaehlerrhosigmak} are fulfiled for the flow $(M, \hat g, \hat\xi)$ with $\hat g=g$ and $\hat\xi=\xi$, we get the following cases:\\
To show $\mathrm{i)}$, we assume that $\mu$ vanishes at one point of $\wit M$. Then, according to \eqref{ONeilltens}, also $h$ vanishes at that point.  Lemma~\ref{l:charactKaehlerrhosigmak}, $\mathrm{ii)}$, then implies that $h$, and thus also $\mu$, vanishes identically. In this case, up to replacing $\rho k$ by $1$ and $\rho\sigma$ by $\rho$, we obtain from the above system that $\rho=\frac{f'}{f'(0)}$.
Furthermore,  up to replacing $u$ by $\frac{1}{f'(0)}u$ which does not change the statement of the theorem, we can assume that $f'(0)=1$. Thus we obtain $\rho=f'$, that is, $\rho(t)=|\wnabla u|\circ F(t,x)$ for all $(t,x)\in\R\times M$ and the triple $(\wit{M}^{2n},\wit{g},\wit{J})$ is locally biholomorphically isometric to $(\R_t\times\R_s\times\Sigma,dt^2\oplus\rho^2(t)ds^2\oplus g_{\Sigma})$, for some K\"ahler manifold $(\Sigma^{2n-2},g_{\Sigma})$ and some positive function $\rho$.  
Note that the eigenvalues $\lambda$ and $\mu$ may be equal (to $0$), which corresponds to the trivial case when $\wnabla u$ is a parallel vector field on $\wit{M}$. \\

To show $\mathrm{ii),(a)}$,   we assume that $\mu$ does not vanish at any point of $\wit M$. 
In this case, it follows from Lemma~\ref{l:charactKaehlerrhosigmak}, $\mathrm{i)}$, that the function $k$ only depends on $t$, so that $(t,x)\mapsto\mu\circ F(t,x)$ only depends on $t$ as well. Setting $\ovl{\rho} :=\rho k$ as well as $\ovl{\sigma}:=\frac{\sigma}{k}$, we can assume that $k\equiv 1$.
Furthermore, Lemma~\ref{l:charactKaehlerrhosigmak}, $\mathrm{iii)}$, implies after rescaling $\xi$ and turning $t$ into $-t$ (which amounts to changing $u$ into $-u$), that $\rho'=\sigma$ and thus $h=-J$ on $M$.
Hence the first equation in the above system allows to get $(\rho\rho')^2=\Big(\frac{f'}{f'(0)}\Big)^2$, that is $\rho\rho'=\frac{f'}{f'(0)}$, or equivalently $\rho^2(t)-\rho^2(0)=\frac{2}{f'(0)}(f(t)-f(0))$, for all $t\in\R$.
Thus, up to replacing $u$ by $\frac{2}{f'(0)}(u-u_0)+\rho^2(0)$, which does not affect the assumptions on $u$, we may assume that $f'(0)=2$ and $f(0)=\rho(0)^2=1$; in particular we obtain $\rho^2=f$ and thus $\rho=\sqrt{f}=\sqrt{u\circ F}$. \color{black} 
We also deduce that $\lambda\circ F(t,x)=\wnabla^2u(\nu,\nu)\circ F(t,x)=f''(t)$ and that $\mu\circ F(t,x)=\frac{f'(t)\rho'(t)}{\rho(t)}=\frac{(f')^2(t)}{2f(t)}=\frac{|\wnabla u|^2}{2u}\circ F(t,x)$.\\     \color{black}

 Now, we consider the case when $n=2$. When $\mu$ vanishes identically, then as before the triple $(\wit{M}^{2n},\wit{g},\wit{J})$ is locally biholomorphically isometric to $(\R_t\times\R_s\times\Sigma,dt^2\oplus\rho^2(t)ds^2\oplus g_{\Sigma})$, for some surface $(\Sigma^{2},g_{\Sigma})$ and some positive function $\rho$. In the following, we consider the case when $\mu$ does not vanish identically on $\wit M$. \color{black} 
Let $x_0\in\wit M$ be a point where $\mu(x_0)\neq 0$. Up to changing the regular value $u_0$ of $u$, we may assume that $x_0\in M=u^{-1}(\{u_0\})$.
This time, we write $\hat{g}_{\xi^\perp}:=\frac{1}{\beta(x)^2}\cdot g_{\xi^\perp}$ for some positive basic function $\beta$ that will be later determined.  From Lemma \ref{l:transvconfchange}, the triple $(M,\hat{g}:=g_{\xi}\oplus \frac{1}{\beta(x)^2}\cdot g_{\xi^\perp}, \hat\xi:=\xi)$ is still a minimal Riemannian flow and, as $n=2$, is also endowed with the same transversal complex structure $J=\wit J$. 
Now we apply Lemma \ref{l:charactKaehlerrhosigmak} to the flow $(M,\hat{g}, \hat\xi)$ and put $F_t^*\wit{g}$ under the form $dt^2\oplus\rho(t)^2\left(\sigma(t)^2\hat{g}_{\xi}\oplus k(t,x)^2\hat{g}_{\xi^\perp}\right)$, where $\rho,\sigma,k$ satisfy \color{black}
\begin{equation}\label{eq:rhosigmanequal2}
\left\{\begin{array}{ll}\rho^2(t)\sigma^2(t)&=\frac{(f')^2(t)}{(f')^2(0)}\\\rho^2(t) k^2(t,x)&=\beta^2(x)\cdot\exp\left(\displaystyle 2\int_0^t\frac{\mu\circ F(s,x)}{f'(s)}\,ds\right).\end{array}\right.
\end{equation}
Lemma~\ref{l:charactKaehlerrhosigmak}, $\mathrm{i)}$, implies the existence of a basic function $C$ on $M$ such that $\frac{\partial(\rho k)^2}{\partial t}=2\rho\sigma C$, in particular $\rho^2(t)k^2(t,x)-\rho^2(0)k^2(0,x)=2C(x)\cdot\int_0^t\rho(s)\sigma(s)\,ds$ for all $(t,x)\in\R\times M$.
But $\int_0^t\rho(s)\sigma(s)\,ds=\frac{f(t)-f(0)}{f'(0)}$ by the first identity of (\ref{eq:rhosigmanequal2}), so that 
\begin{equation}\label{eq:muC}
\rho^2(t)k^2(t,x)-\rho^2(0)k^2(0,x)=2C(x)\cdot\frac{f(t)-f(0)}{f'(0)}\qquad\forall\,(t,x)\in\R\times M.
\end{equation}
Another consequence of (\ref{eq:rhosigmanequal2}) is that $\mu$ and $C$ have the same sign everywhere: we have the identity $\displaystyle{\int_0^t\frac{\mu\circ F(s,x)}{f'(s)}\,ds=\ln(\rho(t)k(t,x))-\ln(\beta(x))}$, from which we deduce that
\[\mu\circ F(t,x)=f'(t)\cdot\frac{1}{2(\rho k)^2}\cdot\frac{\partial(\rho k)^2}{\partial t}(t,x)=\frac{f'(t)(\rho\sigma)(t)}{(\rho k)(t,x)^2}\cdot C(x)\]
with $\frac{f'(t)(\rho\sigma)(t)}{(\rho k)(t,x)^2}>0$ as we recall that $f'(t)=|\wnabla u|_{F(t,x)}$. \color{black}
Therefore, we deduce that $C(x_0)\neq0$ because of $\mu(x_0)\neq0$.
Note in particular that $f$ (or, equivalently, $u$) is necessarily bounded below or above by identity \eqref{eq:muC}: indeed, as $(\rho(t)k(t,x_0))^2>0$, we get $f(t)\geq f(0)-\frac{(\rho(0)k(0,x_0))^2f'(0)}{2C(x_0)}$, if $C(x_0)>0$ and, if $C(x_0)<0$, we have $f(t)\leq f(0)-\frac{(\rho(0)k(0,x_0))^2 f'(0)}{2C(x_0)}$, for all $t\in\R$. Hence, up to changing $u$ into $-u\pm c$ for some $c\in\R$, we may assume that $u_0=f(0)=\frac{(\rho(0)k(0,x_0))^2 f'(0)}{2C(x_0)}$ and $f>0$. \color{black} 
As in this case $C(x_0)>0$, hence $C>0$ on some nonempty open neighbourhood $U$ of $x_0$ in $M$. Therefore, we may set $\beta(x):=\sqrt{\frac{2C(x)f(0)}{f'(0)}}$, for all $x\in U$. Since $(\rho(0)k(0,x))^2=\beta(x)^2$  from the second identity in \eqref{eq:rhosigmanequal2} evaluated at $t=0$, then  \eqref{eq:muC} implies that
$(\rho(t)k(t,x))^2=2C(x)\frac{f(t)}{f'(0)}$, in particular $k$ is in product form. Again Lemma \ref{l:charactKaehlerrhosigmak}, $\mathrm{iii)}$ together with a possible further rescaling of $\xi$ and a change of $u$ into $\frac{2}{f'(0)}u$ \color{black} yields the splitting result in $\mathrm{ii),(b)}$ on $U$. 
Furthermore, as above, we have $\mu=\frac{|\wnabla u|^2}{2u}$ on $\R\times U$. Now the closed and nonempty subset of $\wit{M}$ where $|\wnabla u|^2=2\mu u$ is open.
Indeed, if $|\wnabla u|^2=2\mu u$ is satisfied at some point $(t,z)\in \R\times M$, then since we know that $u>0$ on $\wit{M}$ we also know that $\mu(z)>0$ and thus $C(z)>0$, because $\mathrm{sgn}(C)=\mathrm{sgn}(\mu)$, and therefore $C>0$ on some open connected neighbourhood $V$ of $z$ in $M$.
Repeating on $V$ the argument performed on $U$, we obtain that $|\wnabla u|^2=2\mu u$ on $\R\times V$. This implies that $|\wnabla u|^2=2\mu u$ holds on all $\wit{M}$ by connectedness.
This concludes the proof of Theorem \ref{t:cardwp2eigenv}.
\findemo

\begin{ecor}\label{c:cartyp1} 
Let $(\wit{M}^{2n},g,J)$ be a complete K\"ahler manifold admitting a function $u\in C^\infty(\wit{M},\R_+^\times)$ with
\beit\item $|\wnabla u|=2u$,
\item $\wnabla^2u=2u\left(\xi^\flat\otimes\xi+\nu^\flat\otimes\nu+\mathrm{Id}_{T\wit{M}}\right)$.
\eeit 
Then $(\wit{M}^{2n},g,J)$ is biholomorphically isometric to $(\R\times M^{2n-1},dt^2\oplus e^{2t}\left(e^{2t}\hat{g}_{\hat{\xi}}\oplus \hat{g}_{\hat{\xi}^\perp}\right))$, where $(M^{2n-1},\hat{g},\hat{\xi})$ is Sasaki.
\end{ecor}

{\it Proof}: Note that, by assumption, $u$ has no critical point (because of the first condition and $u>0$ everywhere), $\wnabla^2u$ has pointwise two eigenvalues, $4u$ and $2u$ with $\ker(\wnabla^2u-4u\mathrm{Id})=\mathrm{Span}(\wnabla u,J\wnabla u)$.
By Theorem \ref{t:cardwp2eigenv} -- that applies since by assumption $\frac{|\wnabla u|^2}{2u}=2u$ -- it suffices to notice that $f(t):=u(F(t,x))=f(0)e^{2t}$.
But this obviously follows from $f'(t)=|\wnabla u|(F(t,x))=2u(F(t,x))=2f(t)$.
Choosing $M:=u^{-1}(\{f(0)\})$, we see that we may choose $f(t)=e^{2t}$ and conclude with Theorem~\ref{t:cardwp2eigenv}.
\findemo



\begin{eexemple}
	{\rm The generalized Taub-NUT metrics of Iwai-Katayama on $\mathbb{C}^2$ as described in \cite[Example 2.2]{MorMor11} are Ricci-flat doubly-warped product K\"ahler metrics and therefore are a particular case of our description in Section \ref{s:Kaehler+Einstein}.}
\end{eexemple}

\section{K\"ahler-Einstein doubly-warped products}\label{s:Kaehler+Einstein}

The purpose of this section is to give a characterization of the K\"ahler doubly-warped products of the form $(\wit{M}=I\times M^{2n-1},\wit{g}:=dt^2\oplus\rho^2((\rho')^2\hat{g}_{\hat{\xi}}\oplus \hat{g}_{\hat{\xi}^\perp}), \widetilde J)$, whose underlying metric $\wit{g}$ is Einstein.  Recall first that $(M,\hat{g},\hat{\xi})$ is a minimal Riemannian flow endowed with a complex structure $J$ and that the complex structure $\wit{J}$ on $\wit{M}$ is always the one given by \eqref{newcomplexstructure}. According to Lemma~\ref{l:charactKaehlerrhosigmak}, $i)$, and since here $k=1$ and $\sigma=\rho'$, the complex structure $\wit{J}$ is K\"ahler on $\wit{M}$ and we have $\hat{h}=-J$; hence $(M,\hat{g},\hat{\xi})$ is a Sasakian manifold. We will show in the sequel that the Einstein condition on $(\wit{M},\wit{g})$ is equivalent to $(M,\hat{g},\hat{\xi})$ being $\eta$-Einstein and $\rho$ satisfying an ODE of order $3$. Depending on the sign of the Einstein constant, we will provide solutions of this ODE that in some cases might not be complete.\\

In the following, we will compute the Ricci curvature of $(\wit{M},\wit{g})$ in terms of the transversal Ricci curvature which is associated to the transversal Levi-Civita connection $\hat{\nabla}$, by using the formulas in Lemma \ref{l:identitiesrhosigmak}. For this, we denote by $\left(e_j\right)_{1\leq j\leq 2n-1}$ a local o.n.b. of $TM$ for the metric $g_t$ with $e_{2n-1}=\xi=\frac{1}{\rho\rho'}\hat{\xi}$.
Then,
with our convention, $\wit{R}_{X,Y}=\left[\wnabla_X,\wnabla_Y\right]-\wnabla_{[X,Y]}$ for all $X,Y$, we compute
\be 
\wit{g}(\wit{R}_{\partial_t,\xi}\xi,\partial_t)&=&-\wit{g}(\wit{R}_{\partial_t,\xi}\partial_t,\xi)=-\wit{g}(\wnabla_{\partial_t}\wnabla_{\xi}\partial_t,\xi)+\wit{g}(\wnabla_{\xi}\underbrace{\wnabla_{\partial_t}\partial_t}_{0},\xi)+\wit{g}(\wnabla_{[\partial_t,\xi]}\partial_t,\xi)\\
&=&-\wit{g}(\wnabla_{\partial_t}(\frac{(\rho\rho')'}{\rho\rho'}\xi),\xi)-\frac{(\rho\rho')'}{\rho\rho'}\wit{g}(\wnabla_\xi\partial_t,\xi)\\
&=&-\left(\frac{(\rho\rho')'}{\rho\rho'}\right)'\underbrace{\wit{g}(\xi,\xi)}_{1}-\frac{(\rho\rho')'}{\rho\rho'}\underbrace{\wit{g}(\wnabla_{\partial_t}\xi,\xi)}_{0}-\left(\frac{(\rho\rho')'}{\rho\rho'}\right)^2\wit{g}(\xi,\xi)\\
&=&-\frac{(\rho\rho')''}{\rho\rho'}
\ee
and, for every $j\in\{1,\ldots,2n-2\}$, we similarly compute: $\wit{g}(\wit{R}_{\partial_t,e_j}e_j,\partial_t)=-\frac{\rho''}{\rho}\wit{g}(e_j,e_j)$.
Therefore 
$$\wit{\mathrm{ric}}(\partial_t,\partial_t)=-(2n-2)\frac{\rho''}{\rho}-\frac{(\rho\rho')''}{\rho\rho'}=-(2n+1)\frac{\rho''}{\rho}-\frac{\rho'''}{\rho'}.$$
Note that, because $(\wit{M},\wit{g},\wit{J})$ is K\"ahler, we also have by $\wit{J}$-invariance of the Ricci-curvature that $\wit{\mathrm{ric}}(\xi,\xi)=\wit{\mathrm{ric}}(\partial_t,\partial_t)$, as well as $\wit{\mathrm{ric}}(\xi,\partial_t)=0$.
For every $Z\in\{\xi,\partial_t\}^\perp$, we now compute $\wit{\mathrm{ric}}(Z,\partial_t)$. Indeed, 
\be
\wit{g}(\wit{R}_{Z,e_j}e_j,\partial_t)&=&-\wit{g}(\wit{R}_{Z,e_j}\partial_t,e_j)=-\wit{g}(\wnabla_Z\wnabla_{e_j}\partial_t,e_j)+\wit{g}(\wnabla_{e_j}\wnabla_Z\partial_t,e_j)+\wit{g}(\wnabla_{[Z,e_j]}\partial_t,e_j)\\
&=&-\frac{\rho'}{\rho}\underbrace{\wit{g}(\wnabla_Ze_j,e_j)}_{0}+\frac{\rho'}{\rho}\wit{g}(\wnabla_{e_j}Z,e_j)+\wit{g}(\wnabla_{[Z,e_j]^Q}\partial_t,e_j)-2\wit{g}([Z,e_j],\xi)\underbrace{\wit{g}(\wnabla_\xi\partial_t,e_j)}_{0}\\
&=&\frac{\rho'}{\rho}\wit{g}([e_j,Z],e_j)+\frac{\rho'}{\rho}\wit{g}([Z,e_j],e_j)=0
\ee
and similarly we obtain $\wit{g}(\wit{R}_{Z,\xi}\xi,\partial_t)=0$,
so that $\wit{\mathrm{ric}}(Z,\partial_t)=0$.
Consequently, by the $\wit{J}$-invariance of Ricci-curvature, we also have $\wit{\mathrm{ric}}(Z,\xi)=0$.
The last term to be computed is $\wit{\mathrm{ric}}(Z,Z)$, for any $Z$. Similarly as above, we find $\wit{g}(\wit{R}_{Z,\partial_t}\partial_t,Z)=\wit{g}(\wit{R}_{Z,\xi}\xi,Z)=-\frac{\rho''}{\rho}\wit{g}(Z,Z)$.
In order to compute the remaining curvature term, we take for simplification $e_j$ and $Z$ to be parallel with respect to the connection $\hat\nabla$ at some point $x$ (in this case $[Z,e_j]_x=-2\hat{g}(\hat{h}Z,e_j)\hat\xi_x$, so it is collinear to $\hat\xi_x$). Then, as $\hat{h}=-J$, we have at the point $x$:
\be
\wit{g}(\wit{R}_{Z,e_j}e_j,Z)&=&\wit{g}(\wnabla_Z\wnabla_{e_j}e_j,Z)-\wit{g}(\wnabla_{e_j}\wnabla_Ze_j,Z)-\wit{g}(\wnabla_{[Z,e_j]}e_j,Z)\\
&=&\wit{g}(\wnabla_Z(\hat\nabla_{e_j}e_j-\frac{\rho'}{\rho}\wit{g}(e_j,e_j)\partial_ t),Z)-\wit{g}(\wnabla_{e_j}(\hat\nabla_Z e_j+\frac{\rho'}{\rho}\wit{g}(JZ,e_j)\xi-\frac{\rho'}{\rho}\wit{g}(Z,e_j)\partial_ t),Z)\\&&-\wit{g}([Z,e_j],\xi)\wit{g}(\wnabla_{\xi}e_j,Z)\\
&=&\wit{g}(\hat\nabla_Z\hat\nabla_{e_j}e_j,Z)-(\frac{\rho'}{\rho})^2\wit{g}(Z,Z)-\wit{g}(\hat\nabla_{e_j}\hat\nabla_Ze_j,Z)-(\frac{\rho'}{\rho})^2\wit{g}(JZ,e_j)^2+(\frac{\rho'}{\rho})^2\wit{g}(Z,e_j)^2\\&&-2(\frac{\rho'}{\rho})^2\wit{g}(Je_j,Z)^2\\
&=&\wit{g}(R_{Z,e_j}^{\hat\nabla} e_j,Z)+\left(\frac{\rho'}{\rho}\right)^2(\wit{g}(Z,e_j)^2-\wit{g}(Z,Z)-3\wit{g}(JZ,e_j)^2).
\ee
For the last identity, recall that $R_{Z,Z'}^{\hat\nabla}=\left[\hat\nabla_Z,\hat\nabla_{Z'}\right]-\hat\nabla_{[Z,Z']}$ for all $Z,Z'\in Q=\{\xi,\partial_t\}^\perp$. Finally, we deduce that 
\be 
\wit{\mathrm{ric}}(Z,Z)&=&\sum_{j=1}^{2n-2}\left(\wit{g}(R_{Z,e_j}^{\hat\nabla} e_j,Z)+\left(\frac{\rho'}{\rho}\right)^2(\wit{g}(Z,e_j)^2-\wit{g}(Z,Z)-3\wit{g}(JZ,e_j)^2)\right)-2\frac{\rho''}{\rho}\wit{g}(Z,Z)\\
&=&\mathrm{ric}^{\hat{\nabla}}(Z,Z)-2\left(\frac{\rho\rho''+n(\rho')^2}{\rho^2}\right)\wit{g}(Z,Z).
\ee
To sum up, the $(1,1)$-Ricci-tensor of $(\wit{M},\wit{g})$ is given pointwise by 
$$\wit{\mathrm{Ric}}=-\left((2n+1)\frac{\rho''}{\rho}+\frac{\rho'''}{\rho'}\right)\cdot\left(dt\otimes\!\partial_t\oplus\xi^\flat\!\otimes\!\xi\right)\oplus\left(\frac{1}{\rho^2}\mathrm{Ric}^{\hat{\nabla}}-2\left(\frac{\rho\rho''+n(\rho')^2}{\rho^2}\right)\cdot\mathrm{Id}_{\{\xi,\partial_t\}^\perp}\right).$$ \color{black}
In particular, the manifold $(I\times M^{2n-1},dt^2\oplus\rho^2((\rho')^2\hat{g}_{\hat{\xi}}\oplus\hat{g}_{\hat{\xi}^\perp}))$ is K\"ahler-Einstein if and only if there exists a constant $C\in\mathbb{R}$, which is equal to $\displaystyle\frac{\wit{\mathrm{Scal}}}{2n}$, such that 
$$\left\{\begin{array}{ll}-(2n+1)\frac{\rho''}{\rho}-\frac{\rho'''}{\rho'}&=C\\ 
\frac{1}{\rho^2}\mathrm{Ric}^{\hat{\nabla}}-2\left(\frac{\rho\rho''+n(\rho')^2}{\rho^2}\right)\cdot\mathrm{Id}_{\{\xi,\partial_t\}^\perp}&=C\cdot  \mathrm{Id}_{\{\xi,\partial_t\}^\perp}      
\end{array}\right.,$$
that is, such that 
\begin{equation}\label{eq:witgEinstein}
\left\{\begin{array}{rl}\rho\rho'''&=-(2n+1)\rho'\rho''-C\rho\rho'\\\mathrm{Ric}^{\hat{\nabla}}&=\left(2(\rho\rho''+n(\rho')^2)+C\rho^2\right)\cdot\mathrm{Id}_{\{\xi,\partial_t\}^\perp}.\end{array}\right.
\end{equation}
Note that, by the first equation, the factor $2(\rho\rho''+n(\rho')^2)+C\rho^2=2c$ is constant on $I$ (its first derivative is twice the difference between left- and right-hand-sides of the first equation). This shows, in particular, that the manifold $(M,\hat{g},\hat{\xi})$ is transversally Einstein and therefore it is $\eta$-Einstein, as it is Sasakian. Notice also that $c$ is an integration constant, which is equal to $\displaystyle\frac{\mathrm{Scal}^Q}{4(n-1)}$. \color{black}\\ 

In order to solve this ODE, we consider the following change of variables $z:= ({\rho'})^2+\varepsilon\rho^2$, where $\varepsilon$ denotes the sign of $\widetilde{\mathrm{Scal}}=2nC$, i.e. it is defined as follows:
\[\varepsilon:=\begin{cases} -1,  &\text{if\, } C<0,\\  \phantom{-}0,  & \text{ if\, } C=0, \\  \phantom{-}1,  &\text{ if\, }  C>0. \end{cases}\]
After rescaling the metric in such a way that $C=2(n+1)\varepsilon$, the derivative $z'=2\rho'\rho''+2\varepsilon\rho\rho'$ can be computed as follows:
\begin{eqnarray*}z'
	&=&2\rho'\left(\frac{c}{\rho}-\frac{C}{2}\rho-n\frac{(\rho')^2}{\rho}+\varepsilon\rho\right)\\
	&=&2\rho'\left(\frac{c}{\rho}-\varepsilon n\rho-n\frac{(\rho')^2}{\rho}\right)\\
	&=&2c\frac{\rho'}{\rho}-2n\frac{\rho'}{\rho}z.
\end{eqnarray*}
Solving this linear first-order ODE in $z$, we obtain 
\[z=\frac{c}{n}+D\rho^{-2n},\]
for some constant $D\in\R$. In turn, this leads to the following nonlinear first-order ODE in $\rho$, since we recall that $\rho'>0$:
\begin{equation}\label{eq:rhoEinsteinorder1}
\rho'=\sqrt{-\varepsilon\rho^2+D\rho^{-2n}+\frac{c}{n}}.
\end{equation}
In the following we assume $\rho$ to be defined at $t=0$ with $\rho(0)>0$.
In order to study the solution of this ODE, we distinguish three cases, according to the sign of the constant $D$, as follows:\\

{\bf 1)} If $D=0$, then the equation (\ref{eq:rhoEinsteinorder1}) can either be solved explicitely or it admits no solution, depending on the sign of $\varepsilon$.
\begin{itemize}
	\item  If $\varepsilon=-1$, i.e. if $\widetilde{\mathrm{Scal}}$ is negative, then the ODE \eqref{eq:rhoEinsteinorder1} becomes $\displaystyle \rho'={\sqrt{\rho^2+\frac{c}{n}}}$ and it always admits an explicit maximal solution, which is given as follows, according to the sign of $c$, i.e. the sign of transversal scalar curvature $\mathrm{Scal}^Q$:
	\begin{itemize}
		\item If $c<0$, then $\displaystyle \rho(t)=\sqrt{-\frac{c}{n}}\cdot\cosh\left(t+\mathrm{argcosh}\left(\rho(0)\sqrt{-\frac{n}{c}}\right)\right)$ solves the ODE \eqref{eq:rhoEinsteinorder1} on the interval $I_{\rm{max}}=\displaystyle\left]-\mathrm{argcosh}\left(\rho(0)\sqrt{-\frac{n}{c}}\right),\infty\right[$. For establishing this maximal interval of definition, we use the fact that both functions $\rho$ and $\rho'$ must be positive.
		\item If $c=0$, then $\rho(t)=\rho(0)\cdot e^t$ solves the ODE \eqref{eq:rhoEinsteinorder1} on $\R$.
		\item If $c>0$, then $\displaystyle \rho(t)=\sqrt{\frac{c}{n}}\,\sinh\left(t+\mathrm{argsinh}\left(\rho(0)\sqrt{\frac{n}{c}}\right)\right)$ solves the ODE \eqref{eq:rhoEinsteinorder1} on the interval $I_{\rm{max}}=\displaystyle \left]-\mathrm{argsinh}\left(\rho(0)\sqrt{\frac{n}{c}}\right),\infty\right[$.
	\end{itemize}	
	\item If $\varepsilon=0$, i.e. if $\widetilde{\mathrm{Scal}}$ vanishes,  then the ODE \eqref{eq:rhoEinsteinorder1} becomes $\displaystyle \rho'=\sqrt{\frac{c}{n}}$. Hence it has no solution if $c\leq 0$. However for $c>0$, the function $\displaystyle \rho(t)=t\sqrt{\frac{c}{n}}+\rho(0)$ solves the ODE on $\R$, but $\rho$ is positive only on the interval $I_{\rm{max}}=\left]-\sqrt{\frac{n}{c}}\rho(0), \infty\right[$.
	\item If $\varepsilon=1$, i.e. if $\widetilde{\mathrm{Scal}}$ is positive, then the ODE \eqref{eq:rhoEinsteinorder1} becomes $\displaystyle \rho'=\sqrt{\frac{c}{n}-\rho^2}$. Hence it has no solution with positive derivative, if $c\leq 0$. For $c>0$, the function $\displaystyle \rho(t)=\sqrt{\frac{c}{n}}\,\sin\left(t+\mathrm{arcsin}\left(\rho(0)\sqrt{\frac{n}{c}}\right)\right)$ solves the equation on the interval $I_{\rm{max}}=\displaystyle \left]-\mathrm{arcsin}\left(\rho(0)\sqrt{\frac{n}{c}}\right),\frac{\pi}{2}-\mathrm{arcsin}\left(\rho(0)\sqrt{\frac{n}{c}}\right)\right[$.
\end{itemize}
\smallskip

{\bf 2)} If $D<0$, then we consider the function 
$$f\colon]0,\infty[\to\R, \; f(x):=-\varepsilon x^2+Dx^{-2n}+\frac{c}{n},$$
whose derivative is given by $f'(x)=-2x(\varepsilon+nDx^{-2n-2})$, for all $x>0$. According again to the sign of $\varepsilon$, we distinguish the following three subcases:
\begin{itemize}
	\item If $\varepsilon=-1$, then $f'(x)>0$, for all $x>0$, so the function $f$ is  increasing on $]0,\infty[$ with $\displaystyle \lim_{x\to 0^+}f(x)=-\infty$ and  $\displaystyle\lim_{x\to \infty}f(x)=\infty$. Hence, there exists a unique $\rho_0\in\,]0,\infty[$ with $f(x)<0$ for $0<x<\rho_0$, $f(\rho_0)=0$ and $f(x)>0$ for $x>\rho_0$. Necessarily the solution $\rho$ satisfies $\rho\geq\rho_0$, and actually $\rho>\rho_0$ unless $\rho$ is constant, which is excluded because $\rho'$ is positive everywhere.
	Integrating from $0$ to some positive $t$, we obtain
	\[\int_{\rho(0)}^{\rho(t)}\frac{d\rho}{\sqrt{\rho^2+D\rho^{-2n}+\frac{c}{n}}}=t.\]
	Since $\displaystyle \frac{1}{\sqrt{f(x)}}\mathop\sim_{\rho_0}\frac{1}{\sqrt{f'(\rho_0)\cdot(x-\rho_0)}}$, the solution $\rho$ runs backward to $\rho_0$ in finite time. On the other hand, because of $\displaystyle \frac{1}{\sqrt{f(x)}}\mathop\sim_{\infty}\frac{1}{x}$, the solution $\rho$ becomes infinite only in infinite time. 
	Summing up, we conclude that in this case the maximal solution $\rho$ of the ODE \eqref{eq:rhoEinsteinorder1} is defined on some time interval of the form $I_{\rm{max}}=]t_0,\infty[$,  where $t_0<0$, and fulfils $\displaystyle \lim_{t\to t_0^+}\rho(t)=\rho_0$, $\displaystyle\lim_{t\to \infty}\rho(t)=\infty$ and $\displaystyle \lim_{t\to t_0^+}\rho'(t)=0$, because $f(\rho_0)=0$.
	\item If $\varepsilon=0$, then $f'(x)=-2nDx^{-2n-1}>0$, for all $x>0$, so the function $f$ is increasing with $\displaystyle \lim_{x\to 0^+}f(x)=-\infty$ and  $\displaystyle\lim_{x\to \infty}f(x)=\frac{c}{n}$. Hence, if $c\leq 0$, there is no solution of the ODE~\eqref{eq:rhoEinsteinorder1}. If $c>0$, then there exists a unique $\rho_0$ with $f(x)<0$ for $0<x<\rho_0$, $f(\rho_0)=0$ and $f(x)>0$ for $x>\rho_0$. The same argument as in the previous case shows the solution $\rho$ necessarily satisfies $\rho>\rho_0$. Integrating again from $0$ to some positive $t$, we obtain
	\[\int_{\rho(0)}^{\rho(t)}\frac{d\rho}{\sqrt{D\rho^{-2n}+\frac{c}{n}}}=t.\]
	Since $\displaystyle \frac{1}{\sqrt{f(x)}}\mathop\sim_{\rho_0}\frac{1}{\sqrt{f'(\rho_0)\cdot(x-\rho_0)}}$, the solution $\rho$ runs backward to $\rho_0$ in finite time. Summing up, we conclude that in this case the maximal solution $\rho$ of the ODE \eqref{eq:rhoEinsteinorder1} is defined on some time interval of the form $]t_0,\infty[$, where $t_0<0$, and fulfils $\displaystyle \lim_{t\to t_0^+}\rho(t)=\rho_0$, $\displaystyle\lim_{t\to \infty}\rho(t)=\infty$ and $\displaystyle \lim_{t\to t_0^+}\rho'(t)=0$.
	\item If $\varepsilon=1$, then the function $f$ defined above attains its maximum at $\rho_0=(-nD)^{\frac{1}{2n+2}}$, where the function takes the value $f(\rho_0)=\frac{c}{n}-\frac{n+1}{n}(-nD)^{\frac{1}{n+1}}$. This leads us to consider the following subcases, according to the sign of $f(\rho_0)$, or equivalently, according to the value of $c$, as follows:
	\begin{itemize}
		\item If $c\leq(n+1)(-nD)^{\frac{1}{n+1}}$, then $f(x)\leq f(\rho_0)\leq 0$, hence in this case the ODE~\eqref{eq:rhoEinsteinorder1} has no solution with positive derivative.
		\item  If $c>(n+1)(-nD)^{\frac{1}{n+1}}$, then there exist $\rho_1,\rho_2$, such that $0<\rho_1<\rho_0<\rho_2$ and $f(\rho_1)=f(\rho_2)=0$, which implies that in this case the solution $\rho$ of the ODE~\eqref{eq:rhoEinsteinorder1} is bounded as follows: $\rho_1< \rho(t)<\rho_2$. Integrating again from $0$ to some positive $t$, we obtain
		\[\int_{\rho(0)}^{\rho(t)}\frac{d\rho}{\sqrt{-x^2+D\rho^{-2n}+\frac{c}{n}}}=t.\]
		Since $\displaystyle \frac{1}{\sqrt{f(x)}}\mathop\sim_{\rho_1}\frac{1}{\sqrt{f'(\rho_1)\cdot(x-\rho_1)}}$, the solution $\rho$ runs backward to $\rho_1$ in finite time. Similarly, since $\displaystyle \frac{1}{\sqrt{f(x)}}\mathop\sim_{\rho_2}\frac{1}{\sqrt{f'(\rho_2)\cdot(x-\rho_2)}}$, the solution $\rho$ runs to $\rho_2$ in finite time. Summing up, we conclude that in this case the maximal solution $\rho$ of the ODE \eqref{eq:rhoEinsteinorder1} is defined on some time interval of the form $I_{\rm{max}}=]t_1,t_2[$, where $t_1<0<t_2$, and fulfils $\displaystyle \lim_{t\to t_1^+}\rho(t)=\rho_1$ and $\displaystyle \lim_{t\to t_2^-}\rho(t)=\rho_2$ .
	\end{itemize}
\end{itemize}

{\bf 3)} If $D>0$, then one may proceed similarly to the analysis in the case { 2)}, also by remarking that the function $-f(x)=\varepsilon x^2-Dx^{-2n}-\frac{c}{n}$ is the same as in case  {2)}, as $-D<0$.

Summing up, we have shown the following result:

\begin{ethm}\label{t:caracdwpEinstein}
	Let $(\wit{M}^{2n},g)=(I\times M^{2n-1},dt^2\oplus\rho^2((\rho')^2g_{\hat{\xi}}\oplus g_{\hat{\xi}^\perp}))$ be a K\"ahler doubly-warped product, where $I$ is an open interval containing $0$, $\rho,\rho'\colon I\to\R$ are positive functions and $(M,\hat{g},\hat{\xi})$ is Sasaki.
	Then the following assertions hold:
	\begin{enumerate}
		\item The manifold $(\wit{M}^{2n},g)$ is an Einstein manifold with Einstein constant $2(n+1)\varepsilon$, where $\varepsilon\in\{-1,0,1\}$, if and only if there exist constants $c,D\in\R$, such that
		\[\left\{\begin{array}{ll}\rho'&=\sqrt{-\varepsilon\rho^2+D\rho^{-2n}+\frac{c}{n}}\\\mathrm{Ric}^{\hat{\nabla}}&=2c\cdot\mathrm{Id}_Q,\end{array}\right.\]
		where $Q:=\{\xi,\nu\}^\perp$ and $\hat{\nabla}$ denotes the transverse connection on $Q$.
		\item In each of the following two cases \color{black}: $D=0,\varepsilon=-1$ and $c=0$, or $D>0,\varepsilon=-1$, $c=-(n+1)(nD)^{\frac{1}{n+1}}$ and $\rho(0)>(nD)^{\frac{1}{2n+2}}$, there exists a solution $\rho$ of the ODE $\rho'=\sqrt{\rho^2+D\rho^{-2n}+\frac{c}{n}}$ which is defined on $\R$. For any other values of $c$, $D$ and $\varepsilon$, there exists a solution $\rho$ of $\rho'=\sqrt{\rho^2-D\rho^{-2n}+\frac{c}{n}}$ defined on a maximal interval $I_{\max}\subsetneq\R$ around $0$.
	\end{enumerate}
\end{ethm}

\appendix
\section{Appendix}
Let us recall here some general facts on warped product structures induced by smooth functions.
The local version of the following can be found in the beautiful paper \cite{Kuehnel}, see in particular \cite[Sec. D]{Kuehnel}.
In the following proposition, the Levi-Civita-connection of $(M^n,g)$ is denoted by $\nabla$ (and so differently from the K\"ahler setting, where $\nabla$ denotes the natural connection on the transverse distribution on $(M,g,\xi)$).

\begin{prop}\label{p:warpedprodu}
Let $(M^n,g)$ be a connected complete Riemannian manifold.
Assume that some $u\in C^\infty(M,\R)$ has no critical point on $M^n$ and satisfies $\nabla^2u(\nabla u)=\lambda\nabla u$ for some $\lambda\in C^\infty(M,\R)$.
Then the manifold $(M^n,g)$ is isometric to $(\R\times\Sigma,dt^2\oplus g_t)$, where $\Sigma$ is a level hypersurface of $u$ and $(g_t)_{t\in\R}$ is a one-parameter-family of Riemannian metrics on $\Sigma$.
\end{prop}

{\it Proof}: Fix $u_0\in u(M)$ and let $\Sigma:=u^{-1}(\{u_0\})\subset M$.
By assumption, $\Sigma$ is a smooth hypersurface in $M$.
Consider the map $F:\R\times\Sigma\to M$ given by the flow of $\nu$, i.e., $f(t,x):=F_t^\nu(x)$ for all $(t,x)\in\R\times\Sigma$.
Note that the flow $(F_t)_t=(F_t^\nu)_t$ is well-defined on $\R$ since $\nu$ is a bounded vector field on the complete Riemannian manifold $(M^n,g)$.
We show that $F$ provides the desired isometry.\\
First, $F$ is a local diffeomorphism: for any $(t,x)\in\R\times\Sigma$ and $(T,X)\in \R\times T_x\Sigma$, one has 
\be 
d_{(t,x)}F(T,X)&=&T\frac{\partial F}{\partial t}(t,x)+d_xF_t^\nu(X)\\
&=&T\nu_{F(t,x)}+d_xF_t^\nu(X)\\
&=&Td_xF_t^\nu(\nu_x)+d_xF_t^\nu(X)\textrm{ since }(F_t^\nu)_*\nu=\nu\\
&=&d_xF_t^\nu(T\nu_x+X),
\ee
so that $d_{(t,x)}F(T,X)=0$ iff $T=0$ and $X=0$ (for $F_t^\nu:M\to M$ is a diffeomorphism).
This shows the invertibility of $d_{(t,x)}F$ and hence that $F$ is a local diffeomorphism.\\
In particular $F(\R\times\Sigma)$ is open in $M$.
But this also implies that $F(\R\times\Sigma)$ is closed in $M$: for one may define the equivalence relation $\sim$ on $M$ via: $x,y\in M$, $x\sim y$ iff there exists a $\hat{u}\in u(M)$ such that $x,y\in F(\R\times u^{-1}(\{\hat{u}\}))$, where $F$ is defined by the flow of $\nu$ (starting this time from the hypersurface $u^{-1}(\{\hat{u}\})$ of $M$).
By the preceding argument, each equivalence class is open in $M$ and hence also closed in $M$.
Since $M$ is connected, this yields $F(\R\times\Sigma)=M$, i.e., $F$ is surjective.
The injectivity of $F$ follows easily from the fact that, for any $x\in\Sigma$, the function $f_x:=u\circ F(\cdot,x):\R\to\R$ is monotonously increasing, for it is smooth with $f_x'(t)=|\nabla u|_{F(t,x)}>0$ for all $t\in\R$: if $F(t,x)=F(t',x')$ for some $(t,x),(t',x')\in\R\times\Sigma$, then the point $F(t,x)$ and $F(t',x')$ lie on the same integral curve of $\nu$; but by the injectivity of $f_x$, that curve intersects $\Sigma$ only in $x$, hence $x=x'$ and, again by the injectivity of $f_x$, it follows $t=t'$.
On the whole, $F$ is a diffeomorphism.
In particular, $\Sigma$ itself must be connected.\\
We now look at the pull-back metric $F^*g$ on $\R\times\Sigma$.
Obviously, $(F^*g)(\frac{\partial}{\partial t},\frac{\partial}{\partial t})=1$ since $\nu$ is a unit vector field.
Moreover, as noticed in \cite[Prop. 2]{Santhanam07}, because $\nabla u$ is a pointwise eigenvector for $\nabla^2u$, the vector field $\nu=\frac{\nabla u}{|\nabla u|}$ is geodesic.
This has the important consequence for the splitting of the metric: for any $(t,x)\in\R\times\Sigma$ and $X\in T_x\Sigma$, we have $(F^*g)_{(t,x)}(\frac{\partial}{\partial t},X)=g_{F(t,x)}(\nu_{F(t,x)},d_xF_t(X))=g_{F(t,x)}(d_xF_t^\nu(\nu_x),d_xF_t(X))=(F_t^\nu)^*g(\nu_x,X)$, where
$$ 
\frac{\partial}{\partial s}(F_s^*g)(\nu_x,X)_{|_{s=t}}=(\mathcal{L}_{\nu}g)((F_t)_*\nu,(F_t)_*X)_{F(t,x)}=(\mathcal{L}_{\nu}g)(\nu,(F_t)_*X)_{F(t,x)}.
$$
But for all $X\in TM$,
\[(\mathcal{L}_{\nu}g)(\nu,X)=g(\nabla_\nu\nu,X)+g(\nabla_X\nu,\nu)=0\]
by the fact that $\nu$ is geodesic of constant length.
Therefore, $\frac{\partial}{\partial s}(F_s^*g)(\nu_x,X)_{|_{s=t}}=0$ for all $t\in\R$, thus $(F_t^*g)(\nu_x,X)=(F_0^*g)(\nu_x,X)=g(\nu_x,X)=0$ for all $(t,x)\in\R\times\Sigma$.
This proves the splitting $F^*g=dt^2\oplus g_t$, where $g_t:=(F_t)^*g_{|_{T\Sigma\times T\Sigma}}$. We note an important consequence of the splitting $F^*g=dt^2\oplus g_t$, namely that the flow $(F_t^\nu)_t$ preserves the level hypersurfaces of $u$, or equivalently, that the function $f_x=u\circ F(\cdot,x)$ defined above actually does not depend on $x$.
For given any further $y\in\Sigma$, consider any smooth curve $c:[0,1]\to\Sigma$ with $c(0)=x$ and $c(1)=y$.
For a fixed $t\in\R$, look at the smooth function $h(s):=u\circ F(t,c(s))$, $s\in [0,1]$.
Its first derivative is given by 
\be 
h'(s)&=&g_{F(t,c(s))}(\nabla u,d_{c(s)}F_t(\dot{c}(s)))\\
&=&|\nabla u|_{F(t,c(s))}g_{F(t,c(s))}(\nu_{F(t,c(s))},d_{c(s)}F_t(\dot{c}(s)))\\
&=&|\nabla u|_{F(t,c(s))}(F_t^*g)(\nu_{c(s)},\dot{c}(s))\\
&=&0,
\ee
so that $h$ is constant and hence $h(0)=u(F(t,x))=h(1)=u(F(t,y))$, i.e., $f_x(t)=f_y(t)$. This concludes the proof.
\findemo


{\bf Acknowledgment:} This long-term project benefited from the generous support of the Universities of Stuttgart, Lorraine, Regensburg and from the conference \emph{Riemann and K\"ahler geometry} held at IMAR (Bucharest) between April 15-19, 2019.
The second named author would like to thank the Alexander von Humboldt foundation and the DAAD for financial support. We would like to thank especially Sergiu Moroianu and the GDRI Eco-Math for his support and interest in our project and for introducing us the Taub-NUT example.
Finally, we are very grateful to  Paul-Andi Nagy for pointing  out to us the relation of our work to other 
classification results on K\"ahler manifolds.

\providecommand{\bysame}{\leavevmode\hbox to3em{\hrulefill}\thinspace}


\begin{thebibliography}{10}
\addcontentsline{toc}{chapter}{References}


\bibitem{ACG}
V.~Apostolov, D.M.J.~Calderbank and P.~Gauduchon, 
\emph{Hamiltonian 2-forms in K\"ahler geometry. I. General theory}, 
J. Differential Geom. \textbf{73} (2006), no. 3, 359--412. 

\bibitem{ACGT}
V.~Apostolov, D.M.J.~Calderbank and P.~Gauduchon, C.T\o nnesen-Friedman, 
\emph{Hamiltonian 2-forms in K\"ahler geometry. II. Global classification}, 
J. Differential Geom. \textbf{68} (2004), no. 2, 277-345. 

\bibitem{Baier97} P.D.~Baier, \emph{\"Uber den Diracoperator auf Mannigfaltigkeiten  mit Zylinderenden}, Diplomarbeit, Universit\"at Freiburg, 1997.

\bibitem{BishopONeill69} R.L.~Bishop, B.~O'Neill, \emph{Manifolds of negative curvature}, Trans. Amer. Math. Soc. \textbf{145} (1969), 1--49.

\bibitem{Carriere84}
Y.Carri\`ere, \emph{Structure transverse des feuilletages}, Toulouse, Ast\'erique \textbf{116} (1984), 31--52.


\bibitem{ChiossiNagy10}
S.G.~Chiossi and P.-A.~Nagy, \emph{Complex homothetic foliations on K\"ahler manifolds}, Bull. Lond. Math. Soc. \textbf{44} (2012), no. 1, 113--124.

\bibitem{DM}
A.~Derdzinski and  G.~Maschler,
\emph{Special K\"ahler-Ricci potentials on compact K\"ahler manifolds},
J. Reine Angew. Math. \textbf{593} (2006), 73--116.

\bibitem{GinHabPilSemm19}
N.~Ginoux, G.~Habib, M.~Pilca and  U.~Semmelmann, \emph{An Obata-type characterization of Calabi metrics on line bundles}, preprint.

\bibitem{GinSemm11}
N.~Ginoux and U.~Semmelmann, \emph{Imaginary K\"ahlerian Killing spinors I}, Ann. Glob. Anal. Geom. \textbf{40} (2011), no. 4, 467--495.








\bibitem{Kuehnel}
W.~K\"uhnel, \emph{Conformal transformations between Einstein spaces}, Conformal geometry (Bonn, 1985/1986), 105--146, Aspects Math., E12, Vieweg, Braunschweig, 1988.

\bibitem{MMP}
F.~Madani, A.~Moroianu and M.~Pilca,
\emph{Conformally related K\"ahler metrics and the holonomy of lcK manifolds},
J. Eur. Math. Soc. \textbf{22} (2020) (1), 119--149.


\bibitem{MolzonPinneyMortensen93}
R.~Molzon and K.~Pinney Mortensen, \emph{A characterization of complex projective space up to biholomorphic isometry}, J. Geom. Anal. \textbf{7} (1997), no. 4, 611--621. 

\bibitem{MorMor11}
A.~Moroianu and S.~Moroianu,\emph{The Dirac operator on generalized Taub-NUT spaces}, Comm. Math. Phys. \textbf{305} (2011), no. 3, 641--656.


\bibitem{Oneil66}
B.~O'Neill,\emph{The fundamental equations of a submersion}, Mich. Math. J. \textbf{305} (1966), 13, 459--469.
\color{black}

\bibitem{Obata62}
M.~Obata, \emph{Certain conditions for a Riemanian manifold to be isometric
with a sphere}, J. Math. Soc. Japan \textbf{14} (1962), 333--340.



\bibitem{RanjSant97}
A.~Ranjan and G.~Santhanam, \emph{A generalization of Obata's theorem}, J. Geom. Anal. \textbf{7} (1997), no. 3, 357--375.

\bibitem{Santhanam07}
G.~Santhanam, \emph{Obata's theorem for K\"ahler manifolds},  Illinois J. Math.  \textbf{51} (2007),  no. 4, 1349--1362. 

\bibitem{Tondeur88}
Ph.~Tondeur, \emph{Foliations on Riemannian manifolds},  Springer, New York, 1988.
\end{thebibliography}
\end{document}